\documentclass[12pt, reqno]{amsart}
\usepackage{amsmath, amstext, amsbsy, amssymb, amscd}
\usepackage{amsmath}
\usepackage{amsxtra}
\usepackage{amscd}
\usepackage{amsthm}
\usepackage{amsfonts}
\usepackage{amssymb}
\usepackage{eucal}
\usepackage{color}

\newtheorem{theorem}{Theorem}[section]
\newtheorem{lemma}{Lemma}[section]
\newtheorem{definition}{Definition}[section]
\newtheorem{prop}{Proposition}[section]
\newtheorem{cor}{Corollary}[section]

\theoremstyle{remark}
\newtheorem{remark}{Remark}[section]

\newcommand{\binf}{{b}_\infty}
\newcommand{\C}{\mathbb C}
\newcommand{\cinf}{ c_{\infty} }
\newcommand{\dinf}{ d_{\infty} }
\newcommand{\gl}{ {\mathfrak {gl}} }
\newcommand{\hgl}{ \widehat{ \mathfrak{gl} } }
\newcommand{\glpm}{ \widehat{{ \mathfrak {gl} }}_{\pm} }
\newcommand{\hf}{ \frac12}
\newcommand{\half}{ \hf}
\newcommand{\hL}{{\Lambda} }
 \newcommand{\vac}{ |0 \rangle }

\newcommand{\la}{\lambda}
\newcommand{\Z}{\mathbb Z}
\newcommand{\HZ}{\frac{1}{2}+\Z}
\newcommand{\fock}{\mathcal F}
\newcommand{\trace}{\,{\rm tr}\,}

\newcommand{\PA} {\mathcal P }
\newcommand{\OSP} {\mathcal{OSP}}
\newcommand{\st} {\,|\,}
\newcommand{\ep}{\varepsilon}
\newcommand{\eps}{\epsilon}
\newcommand{\no} {\text{:}}
\newcommand{\bff}{\mathbf{F}}
\newcommand{\vphi}{\varphi}
\newcommand{\A}{\mathsf A}
\newcommand{\D}{\mathsf D}
\newcommand{\B}{\mathsf B}
\newcommand{\CC}{\mathsf C}

\begin{document}

\title[Bloch-Okounkov functions of classical type]
{The Bloch-Okounkov correlation functions of classical type}

\author{David G. Taylor}
\address{Department of Mathematics, University of Virginia, Charlottesville, VA 22904}
\email{dgtaylor@virginia.edu}

\author[Weiqiang Wang]{Weiqiang Wang} \email{ww9c@virginia.edu}

\begin{abstract}
Bloch and Okounkov introduced an $n$-point correlation function on
the infinite wedge space and found an elegant closed formula in
terms of theta functions. This function has connections to
Gromov-Witten theory, Hilbert schemes, symmetric groups, etc, and it
can also be interpreted as correlation functions on integrable
$\hgl_\infty$-modules of level one. Such $\hgl_\infty$-correlation
functions at higher levels were then calculated by Cheng and Wang.

In this paper, generalizing the type $A$ results, we formulate and
determine the $n$-point correlation functions in the sense of
Bloch-Okounkov on integrable modules over classical Lie subalgebras
of $\hgl_\infty$ of type $B,C,D$ at arbitrary levels. As byproducts,
we obtain new $q$-dimension formulas for integrable modules of type
$B,C,D$ and some fermionic type $q$-identities.
\end{abstract}

\maketitle

\tableofcontents

\section{Introduction}
\subsection{The earlier works}

Bloch and Okounkov \cite{BO} (also see \cite{Ok} for some
simplification) introduced an $n$-point correlation function on the
infinite wedge space and found an elegant closed formula in terms of
theta functions. Their work was in part motivated by certain modular
invariance property of trace functions of vertex operators and
representation theory of the $W_{1+\infty}$ algebra (cf. \cite{Zhu,
FKRW, Bl}). Subsequently, this function and its variant have been
interpreted as a generating function of the Gromov-Witten invariants
of an elliptic curve by Okounkov-Pandharipande \cite{OP}, and  as a
generating function of intersection numbers on Hilbert schemes of
points by Li, Qin and the second author \cite{LQW}. We also refer
the reader to \cite{Lep, Mil} for formal vertex operator
generalizations, \cite{W2} for a neutral fermionic Fock space
version, and \cite{CW2} for a $q,t$-deformation of the
Bloch-Okounkov $n$-point function.

From a representation theoretic viewpoint, the Bloch-Okounkov
$n$-point function can be also easily interpreted as correlation
functions on integrable modules over Lie algebra $\hgl_\infty$ of
level one (cf. \cite{Ok, Mil, CW}). Along this line, Cheng and the
second author \cite{CW} formulated and calculated such $n$-point
correlation functions on integrable $\hgl_\infty$-modules of level
$l$ ($l \in \mathbb N$).
%
%
%
%
\subsection{The goal}
The goal of this paper is to formulate and determine the $n$-point
correlation functions in the sense of Bloch-Okounkov on integrable
modules over classical Lie subalgebras of $\hgl_\infty$ of type
$B,C,D$ at arbitrary level, generalizing the works \cite{BO, CW}
in type $A$. Note that the integrability of these modules implies
that the levels have to be positive (half-)integral. By the
original works of Date, Jimbo, Kashiwara and Miwa (\cite{DJKM1,
DJKM2}), $\hgl_\infty$ affords classical Lie subalgebras of type
$B, C, D$, and these infinite-dimensional Lie algebras played an
important role in connections with solition equations discovered
by the Kyoto school in early 1980's.

The representation theory of $\hgl_\infty$ is intimately related to
that of the $W_{1+\infty}$ algebra (cf. \cite{FKRW} and the
references therein). It follows that the Bloch-Okounkov correlation
functions for $\hgl_\infty$-modules can be regarded as those for
$W_{1+\infty}$-modules. In the same vein, the representation theory
of the classical Lie subalgebras of $\hgl_\infty$ is closely related
to that of the classical Lie subalgebras of $W_{1+\infty}$ initiated
in \cite{KWY}; the Howe dualities \cite{W1}, which are to be used
extensively in this work, readily carry over if one replaces
classical Lie subalgebras of $\hgl_\infty$ by classical Lie
subalgebras of $W_{1+\infty}$. In this way, the $n$-point
correlation functions studied in this paper can be in turn regarded
as those for modules over classical Lie subalgebras of
$W_{1+\infty}$.
\subsection{Our approach}

To achieve our goal, the first (main) step here is to relate the
correlation functions at higher levels to the correlation
functions at the bottom levels (i.e. of level one and/or level
$\half$). Our main tool is the free field realization \cite{DJKM1,
DJKM2} (also cf. Feingold-Frenkel \cite{FF}) and the Howe duality
due to the second author \cite{W1} between the classical Lie
subalgebras of $\hgl_\infty$ and various classical Lie groups
(where sometimes disconnected groups and different covering groups
are required). We refer to \cite{H1, H2} for Howe's original
setups, where all Lie algebras and groups involved are
finite-dimensional.

We note that all integrable modules of these Lie subalgebras of
$\hgl_\infty$ appear in these Howe duality decompositions, and the
level of an integrable module matches with the rank of the
corresponding Lie group. A detailed knowledge of irreducible
modules over various Lie groups (cf. Br\"ocker-tom Dieck [BtD])
and the determinantal ratio form of the Weyl character formulas
for classical Lie algebras (cf. Fulton-Harris \cite{FH}) are also
used in this paper in an essential way.

A similar approach has actually been applied in \cite{CW}
successfully where the type $A$ Howe duality between $\hgl_\infty$
and $GL_l$ due to I.~Frenkel \cite{Fr} (also cf. \cite{W1}) was
used. Forced by the new technical features in type $B,C,D$, we
establish in this paper the relations between the $n$-point
functions at higher levels and at the bottom levels in a different
way, avoiding the usage of the $q$-dimension formula for
integrable modules in \cite{CW}. As a byproduct, we obtain neat
$q$-dimension formulas for the corresponding integrable modules
over the classical Lie subalgebras of type $B,C,D$, which are
simpler than the ones in \cite{KWY} obtained by a specialization
of the Weyl-Kac character formulas. We remark that the idea of
using Howe duality to obtain irreducible character formulas has
also been applicable in different setups (cf. Cheng-Lam
\cite{CL}).

Our second step is more straightforward. By using the free field
realization we are able to relate the calculation of the $n$-point
function of type $B, D$ of level one to the $n$-point function of
type $A$ of level one which has been computed in \cite{BO}. The
type $C$ level one case can be handled by a combination of Howe
duality and the connection to the type $A$ level one case.

An additional step is needed to take care of the half-integral
levels, which occur in type $B$ and $D$. Using an identification of
a pair of complex fermions and two neutral fermions, we obtain
formulas, recursive on $n$, of computing the $n$-point functions of
classical type of level $\half$ in terms of those of level one.
Explicit formulas in different forms for the $1$-point functions of
type $B$ and $D$ of level $\half$ were obtained in \cite{W2} using
the method of partition identities. Identifying these different
formulas gives rise to two interesting $q$-identities of fermionic
type. It turns out that these identities have been known with a very
different proof (cf. e.g. \cite{K}).

Combining all these steps together, we have calculated all the
$n$-point correlation functions of classical type. The final
formulas involve the Weyl groups of the Lie groups appearing in
various Howe dualities and the original Bloch-Okounkov function of
type $A$ and level one (which in turn is an expression in terms of
theta functions). Remarkably, the solutions in type $B$ and type
$D$ look almost identical formally though different Lie algebras
and different Howe dualities are involved in different type.
\subsection{Open questions}
The integrable modules whose correlation functions are computed here
are occasionally not irreducible (instead it could be a sum of two
irreducibles) over the infinite-dimensional Lie algebras of type
$D$, but they can always be regarded as irreducible modules over the
corresponding orthogonal groups. This is a familiar phenomenon of
spinor vs half-spinor representations. Nevertheless, it will be
interesting to determine completely the (refined) $n$-point
functions for all {\em irreducible} integrable modules of type $D$.
In this direction, we have only obtained limited results. By
observing an intrinsic connection with the theory of partitions (cf.
Andrews \cite{Andrews}), we find an explicit formula for the refined
$1$-point function of type $D$ of level one.

We can also formulate the $n$-point correlation functions for
modules of {\em negative} (half-)integral levels of $\hgl_\infty$
and its classical subalgebras, and these modules have appeared in
the Howe duality decompositions (cf. \cite{KR} for type $A$ and
\cite{W1} in general). It will be interesting to determine these
$n$-point correlation functions. In light of the developments in
\cite{CW} and in this paper, the main difficulty lies in
understanding the cases at level $-1$, where the connection with
the theory of partitions available at level one and at level $\hf$
is now lacking.

A more challenging question is to ask for a geometric interpretation
of the correlation functions studied in this paper (and also in
\cite{CW}). For example, can they be interpreted as correlation
functions in some supersymmetric gauge theory where the classical
Lie groups used in this paper appear as gauge groups?
\subsection{Organization and Acknowledgment} The paper is organized as follows.
In Section~\ref{sec:prelim}, we set up the notations for the
classical Lie subalgebras of $\hgl_\infty$ and various Lie groups
used in later sections. In Sections~\ref{sec:dD} and \ref{sec:dB}
respectively, we formulate and calculate the $n$-point functions
and the $q$-dimension of integrable $\dinf$-modules of level $l$
and of level $l+\hf$ respectively. In Section~\ref{sec:cC}, we
calculate the $n$-point functions and the $q$-dimension of
integrable $\cinf$-modules of level $l$. In
Sections~\ref{sec:binfl} and \ref{sec:bB} respectively, we
formulate and calculate the $n$-point functions and the
$q$-dimension of integrable $\binf$-modules of level $l$ and of
level $l+\hf$ respectively.

Various Fock spaces of free fermionic fields and Howe dualities
are recalled and used in each of the
Sections~\ref{sec:dD}--\ref{sec:bB}. The proofs in
Sections~\ref{sec:dD} and \ref{sec:dB} are given in detail, while
the proofs in Sections~\ref{sec:cC}--\ref{sec:bB} are often
sketchy when they are parallel to the ones in
Sections~\ref{sec:dD}--\ref{sec:dB}.

This research is partially supported by NSF and NSA grants. We thank
Shun-Jen Cheng for helpful discussions and comments.

\section{The preliminaries}
\label{sec:prelim}

The purpose of this section is to set up notations for the
infinite-dimensional Lie algebras and classical Lie groups which we
will use.

\subsection{Classical Lie algebras of infinite dimension}
\label{sec_algebras}
In this subsection we review Lie algebras $\hgl \equiv\hgl_\infty$
and its various Lie subalgebras of $B, C, D$ type (cf. \cite{DJKM1,
DJKM2}).
\subsubsection{Lie algebra $\hgl$}

Denote by $\gl$ the Lie algebra of all matrices $(a_{ij})_{i,j \in
\Z}$ satisfying $a_{ij} =0$ for $|i -j|$ sufficiently large.
Denote by $E_{ij}$ the infinite matrix with $1$ at $(i, j)$ place
and $0$ elsewhere and let the weight of $ E_{ij}$ be $j - i $.
This defines a $\Z$--principal gradation $\gl = \bigoplus_{j \in
\Z} \gl_j$. Denote by $\hgl \equiv \hgl_\infty = \gl \oplus \C C$
the central extension given by the following $2$--cocycle with
values in $\C$ (cf. \cite{DJKM1}):
\begin{eqnarray}
 C(A, B) = \trace \left([J, A]B \right)
  \label{eq_cocy}
\end{eqnarray}
where $J = \sum_{ j \leq 0} E_{ii}$. The $\Z$--gradation of Lie
algebra $\gl$ extends to $\hgl$ by letting the weight of $C$ to be
$0$. This leads to a triangular decomposition
$$\hgl = \widehat{\gl}_{+}
          \oplus \widehat{\gl}_{0} \oplus \widehat{\gl}_{-}  $$
where
$ \glpm = \oplus_{ j \in \mathbb N} \widehat{\gl}_{\pm j},
\widehat{\gl}_{0} = \gl_0 \oplus \C C.$
Let
\begin{eqnarray*}
 H^a_i  = E_{ii} - E_{i+1, i+1} + \delta_{i,0} C \quad (i \in \Z).
\end{eqnarray*}
Denote by $L(\hgl; \Lambda)$ the highest weight $\hgl$--module
with highest weight $\Lambda  \in {\hgl}_{0}^{*}$, where $C$ acts
as a scalar which is called the {\em level}. Let $\hL_j^a \in
\hgl_0^{*}$ be the fundamental weights, i.e. $\hL_j^a ( H_i^a ) =
\delta_{ij}.$ The Dynkin diagram for $\hgl$, with fundamental
weights labeled, is the following:
  \begin{equation*}
  \begin{picture}(150,45) 
  \put(-16,23){\dots}
  \put(0,20){$\circ$}
  \put(7,23){\line(1,0){32}}
  \put(40,20){$\circ$}
  \put(47,23){\line(1,0){32}}
  \put(81,20){$\circ$}
  \put(88,23){\line(1,0){32}}
  \put(121,20){$\circ$}
  \put(128,23){\line(1,0){32}}
  \put(161,20){$\circ$}
  \put(170,23){\dots}
  \put(-4,9){$-2$}
  \put(36,9){$-1$}
  \put(81,9){$0$}
  \put(121,9){$1$}
  \put(162,9){$2$}
  \end{picture}
  \end{equation*}
\subsubsection{Lie algebra $\dinf$}
  \label{subsec_dinf}
Let
\begin{eqnarray*}
{\overline{d}}_{\infty} = \{ (a_{ij})_{i,j \in \Z} \in \gl \mid
a_{ij} = -a_{1-j,1-i} \}
 \end{eqnarray*}
be a Lie subalgebra of $\gl$ of type $D$. Denote by $\dinf =
{\overline{d}}_{\infty} \bigoplus \C C $ the central extension
given by the $2$-cocycle (\ref{eq_cocy}). Then $\dinf$ has a
natural triangular decomposition induced from $\hgl$ with Cartan
subalgebra ${\dinf}_0 =  \widehat{\gl}_{0}  \cap \dinf$.  Given
$\Lambda \in {\dinf}_0^* $, we let
\begin{eqnarray*}
  H^d_i & = & E_{ii} + E_{-i, -i} - E_{i+1, i+1} - E_{-i+1, -i+1}
                                \quad (i \in \mathbb N),\\
  H^d_0 & = & E_{0,0} + E_{-1,-1} -E_{2,2} -E_{1,1} + 2C.
\end{eqnarray*}
Denote by $\hL^d_i $ the $i$-th fundamental weight of $\dinf$,
i.e. $\hL^d_i (H^d_j ) = \delta_{ij}$.  The Dynkin diagram of
$\dinf$ is:
  
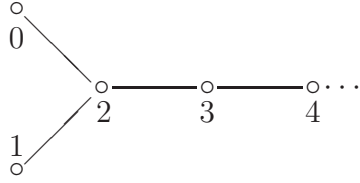
\begin{equation*}
  \begin{picture}(150,75) 
  \put(11,4){\line(1,1){25}}
  \put(5,60){$\circ$}
  \put(11,60){\line(1,-1){25}}
  \put(5,-1){$\circ$}
  \put(37,30){$\circ$}
  \put(44,33){\line(1,0){32}}
  \put(77,30){$\circ$}
  \put(84,33){\line(1,0){32}}
  \put(117,30){$\circ$}
  \put(124,30){$\cdots$}
  \put(5,48){$0$}
  \put(5,7){$1$}
  \put(38,20){$2$}
  \put(77,20){$3$}
  \put(117,20){$4$}
  \end{picture}
  \end{equation*}
\subsubsection{Lie algebra $\cinf$}

Let
\begin{eqnarray*}
 {\overline{c}}_{\infty}
      = \{ (a_{ij})_{i,j \in \Z} \in \gl\mid
                    a_{ij} = - (-1)^{i+j}a_{1-j,1-i} \}
\end{eqnarray*}
be a Lie subalgebra of $\gl$ of type $C$.
Denote by $\cinf$ the central extension of ${\overline{c}}_{\infty}$
given by the $2$-cocycle (\ref{eq_cocy}). Then $\cinf$ inherits from
$\hgl$ a natural triangular decomposition with Cartan subalgebra
${\cinf}_0$. Given $\Lambda \in {\cinf}_0^* $, we let
\begin{eqnarray*}
  H^c_i & = & E_{ii} + E_{-i, -i} - E_{i+1, i+1} - E_{1-i, 1-i}
   \quad (i \in \mathbb N),    \\
   H^c_0 & = & E_{0,0} -E_{1,1} + C.
\end{eqnarray*}
Denote by $\hL^c_i $ the $i$-th fundamental weight of $\cinf$,
i.e. $\hL^c_i (H^c_j ) = \delta_{ij}$. The Dynkin diagram of
$\cinf$ is:
  
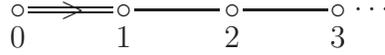
\begin{equation*}
  \begin{picture}(150,45) 
  \put(10,20){$\circ$}
  \put(17,24){\line(1,0){32}}
  \put(17,22){\line(1,0){32}}
  \put(29,20){$>$}
  \put(50,20){$\circ$}
  \put(57,23){\line(1,0){32}}
  \put(91,20){$\circ$}
  \put(98,23){\line(1,0){32}}
  \put(131,20){$\circ$}
  \put(140,23){\dots}
  \put(10,9){$0$}
  \put(50,9){$1$}
  \put(91,9){$2$}
  \put(131,9){$3$}
  \end{picture}
  \end{equation*}
\subsubsection{Lie algebra $\binf$}
   \label{subsec_binf}
 Let
\begin{eqnarray*}
   \overline{b}_{\infty}
   = \{ (a_{ij})_{i,j \in \Z} \in \gl \mid a_{ij} = - a_{-j,-i} \}
\end{eqnarray*}
be a Lie subalgebra of $\gl$ of type $B$. Denote by $\binf$ the
central extension of $\overline{b}_{\infty}$ given by the
$2$-cocycle (\ref{eq_cocy}). The Lie algebra $\binf$ inherits from
$\hgl$ a natural triangular decomposition with Cartan subalgebra
${\binf}_0$. Given $\Lambda \in {\binf}_0^* $, we let
\begin{eqnarray*}
   H^b_i & = & E_{ii} + E_{-i-1, -i-1} - E_{i+1, i+1} - E_{-i, -i}
    \quad (i \in \mathbb N),    \\
   H^b_0 & = & 2 (E_{-1,-1} -E_{1,1}) + 2C.
\end{eqnarray*}
Denote by $\hL^b_i $ the $i$-th fundamental weight of $\binf$,
i.e. $\hL^b_i (H^b_j ) = \delta_{ij}$. The Dynkin diagram of
$\binf$ is:
  
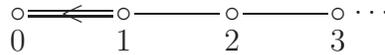
\begin{equation*}
  \begin{picture}(150,45) 
  \put(10,20){$\circ$}
  \put(17,24){\line(1,0){32}}
  \put(17,22){\line(1,0){32}}
  \put(29,20){$<$}
  \put(50,20){$\circ$}
  \put(57,23){\line(1,0){32}}
  \put(91,20){$\circ$}
  \put(98,23){\line(1,0){32}}
  \put(131,20){$\circ$}
  \put(140,23){\dots}
  \put(10,9){$0$}
  \put(50,9){$1$}
  \put(91,9){$2$}
  \put(131,9){$3$}
  \end{picture}
  \end{equation*}
\subsection{Classical Lie groups}
    \label{sec_classical}
We present here a parametrization of irreducible modules of various
classical Lie groups. See \cite{BtD} (also \cite{W1}) for more
detail.

\subsubsection{$O(2l)$}
  \label{subsec_evenorth}

We define $O(2l) = \{ g \in GL(2l); {}^t g J g = J \}$ with
$$J = \begin{bmatrix} 0 & I_l \\
                     I_l & 0 \end{bmatrix}$$
Lie group $GL(l) $ can be identified as a subgroup of $O(2l)$
consisting of matrices of the form $\text{diag}(g,{}^t g^{-1})$,
where ${}^t g$ denotes the transpose. Lie algebra ${\mathfrak
{so}}(2l)$ of $SO(2l)$ consists of matrices of the form
  \begin{eqnarray}
  \left[ \begin{array}{cc}
    \alpha      & \beta       \\
    \gamma      & -{}^t{\alpha}
         \end{array} \right]
  \label{eq_matr}
  \end{eqnarray}
where $\alpha, \beta, \gamma$ are $l \times l$ matrices and $
\beta, \gamma$ are skew-symmetric. Lie algebra ${\mathfrak {gl}}
(l)$ is identified with the subalgebra of ${\mathfrak {so}}(2l)$
consisting of matrices of the form (\ref{eq_matr}) with $\beta
=\gamma =0$. Let ${\mathfrak b} ( {\mathfrak {so}}(2l) )$ be the
Borel subalgebra of ${\mathfrak {so}}(2l)$ which consists of
matrices (\ref{eq_matr}) with $\gamma =0$ and $\alpha$ upper
triangular matrices, and ${\mathfrak h} ({\mathfrak {so}}(2l) )$
be the Cartan subalgebra of diagonal matrices $\text{diag}(t_1,
\ldots, t_l, -t_1, \ldots, - t_l ),$ $t_i \in \C.$ Then
${\mathfrak {gl}}(l)$ and ${\mathfrak {so}}(2l) $ share the same
Cartan subalgebra.

An irreducible $GL(l)$-module is parameterized by its highest
weight which runs over the set
$$\Sigma(A)
\equiv  \{(m_1, m_2, \ldots, m_l ) \mid m_1 \geq m_2 \geq \ldots
\geq m_l ,
 m_i \in \Z \}. $$
An irreducible $SO(2l)$-module is parameterized by its highest
weight in $ \{ (m_1, m_2, \ldots, m_l ),\mid
 m_1 \geq m_2 \geq \ldots \geq  m_{l-1} \geq | m_l |,
 m_i \in \Z \}$. For notational simplicity, we may identify a highest weight
 module with its highest weight below.

$O(2l)$ is a semi-direct product of $SO(2l)$ by $\Z_2$. Denote by
$\tau \in O(2l) \setminus SO(2l)$ the $2l \times 2l$ matrix
\begin{eqnarray}
\left[ \begin{array}{cc}
A & B \\
B & A
\end{array} \right]
  \label{eq_tau}
\end{eqnarray}
with $ A = \text{diag}(1, \ldots, 1, 0), B = \text{diag}(0,
\ldots, 0, 1). $ Then $\tau$ normalizes the Borel $\mathfrak b$.
If $\lambda$ is an $SO(2l)$-module of highest weight $(m_1, m_2,
\ldots, m_l )$, then $\tau . {\lambda}$ has highest weight
$\overline{\la}:=(m_1, m_2, \ldots, - m_l )$. The induced module
of $(m_1, m_2, \ldots, m_l )$ ($m_l \neq 0$) to $O(2l)$ is
irreducible and its restriction to $SO(2l)$ is a sum of $(m_1,
m_2, \ldots, m_l )$ and $(m_1, m_2, \ldots, - m_l )$. We denote
this $O(2l)$-module $\lambda$ by $(m_1, m_2, \ldots, {m}_l )$,
where $m_l > 0$. If $ m_l = 0$, the module $\lambda = (m_1, m_2,
\ldots, m_{l-1}, 0 )$ extends to two different $O(2l)$-modules,
denoted by $\lambda$ and $\lambda \otimes \det$, where $\det$ is
the $1$-dimensional non-trivial $O(2l)$-module. We denote
\begin{eqnarray*}
  \Sigma(D) =
   & & \left\{ (m_1, m_2, \ldots, {m}_l ) \mid
       m_1 \geq m_2 \geq \ldots \geq m_l > 0, m_i \in \Z   \right \}
                               \\
   & & \cup   \left\{ (m_1, m_2, \ldots, m_{l-1}, 0 ) \otimes \det,
   (m_1, m_2, \ldots, m_{l-1}, 0 ) \mid \right.\\
   & &   \left. \quad  m_1 \geq m_2 \geq \ldots
   \geq  m_{l-1} \geq 0, m_i \in \Z \right\}.
\end{eqnarray*}

%

%
%
%
\subsubsection{$O(2l+1)$}

Let $O (2l+1) = \{ g \in GL(2l+1) \mid {}^t g J g = J \}$, where
  \begin{eqnarray*}
  J =\left[ \begin{array}{ccc}
    0      & I_l   & 0 \\
    I_l    & 0     & 0   \\
    0      & 0     & 1
         \end{array} \right].
  \end{eqnarray*}
The Lie algebra ${\mathfrak {so}}(2l+1)$ is the Lie subalgebra of
${\mathfrak {gl}}(2l +1)$ consisting of $(2l+1)\times (2l+1)$
matrices of the form
  \begin{eqnarray}
  \left[ \begin{array}{ccc}
    \alpha      & \beta         & \delta \\
    \gamma      & -{}^t{\alpha} & h      \\
     -{}^t h    & -{}^t{\delta} & 0
         \end{array} \right]
   \label{eq_form}
  \end{eqnarray}
where $\alpha, \beta, \gamma$ are $l \times l$ matrices and $
\beta, \gamma$ skew-symmetric. The Borel subalgebra ${\mathfrak b}
( {\mathfrak {so} } (2l+1) ) $ consists of matrices
(\ref{eq_form}) by putting $\gamma $, $h$, $\delta$ to be $0$ and
$\alpha$ to be upper triangular. The Cartan subalgebra ${\mathfrak
h} ( {\mathfrak {so} } (2l+1) ) $ consists of diagonal matrices of
the form $ \text{diag} ( t_1, \ldots, t_l; -t_1 \ldots -t_l; 0 ),$
$ t_i \in \C$. An irreducible module of $SO(2l+1)$ is
parameterized by its highest weight $(m_1, \ldots, m_l ) \in
\PA^l$, where $\PA^l$ denotes the set of partitions with at most
$l$ non-zero parts.

It is well known that $ O (2l+1)$ is isomorphic to the direct
product $ SO (2l+1) \times \Z_2$ by sending the minus identity
matrix to $-1 \in \Z_2 = \{ \pm 1 \}.$ Denote by $\det$ the
non-trivial one-dimensional representation of $O(2l+1)$. An
representation $\lambda$ of $SO(2l+1)$ extends to two different
representations $\lambda$ and $\lambda \otimes \det$ of $O(2l+1)$.
Then we can parameterize irreducible representations of $O(2l+1)$
by $ (m_1, \ldots, m_l)$ and $ (m_1, \ldots, m_l) \otimes \det $.
We shall denote
$$ \Sigma(B) =
 \PA^l \cup \left\{\la \otimes \det  \mid \la \in \PA^l \right \}.
 $$
\subsubsection{$Spin (n)$ and $Pin (n)$}
  The Pin group $Pin (n)$ is the double covering group
of $ O(n)$, namely we have
$$
1 \longrightarrow \Z_2  \longrightarrow Pin (n)
   \longrightarrow O (n) \longrightarrow 1.
$$
We then define the spin group $Spin (n)$ to be the inverse image
of $SO (n)$ under the projection from $ Pin (n)$ to $ O (n) $.

{\bf Case $ {\bf n = 2l}$.} Let $ {\bf 1}_l = (1, 1, \ldots, 1)$,
${\bar{\bf 1} }_l  = (1, 1, \ldots, 1, -1) \in \Z^l $. An
irreducible representation of $Spin (2l)$ which does not factor to
$SO (2l)$
 is an irreducible representation of ${\mathfrak {so}} (2l)$
parameterized by its highest weight
\begin{eqnarray}
  \lambda = {\bf 1}_l/2 + (m_1, m_2, \ldots, m_l)    \label{eq_wtplus}
\end{eqnarray}
or
\begin{eqnarray}
\overline{\lambda} = {\bar{\bf 1} }_l/2 + (m_1, m_2, \ldots, -m_l)
\label{eq_wtminus}
\end{eqnarray}
where $m_1 \geq \ldots \geq m_l \geq 0, m_i \in \Z.$

There are two possibilities. First, an irreducible representation
of $Pin (2l)$ factors to that of $O(2l)$, then we can use the
parametrization of irreducible representations of $O(2l)$ to
parameterize these representations of $Pin (2l)$.

Secondly, an irreducible representation of $Pin (2l)$ is induced
from an irreducible representation of $Spin (2l)$ with highest
weight of (\ref{eq_wtplus}) or (\ref{eq_wtminus}). When restricted
to $Spin (2l)$, it will decompose into a sum of the two
irreducible representations of highest weights (\ref{eq_wtplus})
and (\ref{eq_wtminus}). We will use $\lambda = {\bf 1}_l  /2 +
(m_1, m_2, \ldots, {m}_l), m_l \geq 0$ to denote this irreducible
representation of $Pin (2l)$. Denote by
$$ \Sigma(Pin) = \{ {\bf 1}_l /2 + (m_1, m_2, \ldots, {m}_l) \mid
   (m_1, m_2, \ldots, {m}_l) \in \PA^l \}.  $$

{\bf Case ${\bf n = 2l +1}$.}
   An irreducible representation of $Spin (2l +1)$
which does not factor to $SO (2l +1)$ is an irreducible
representation of ${\mathfrak {so}} (2l +1)$ parameterized by its
highest weight $\la \in \Sigma(Pin).$
\subsubsection{$Sp (2l)$}

The Lie group $Sp(2l)$ is the subgroup of $GL(2l)$ which preserves
the following skew-symmetric bilinear form
  \begin{eqnarray*}
  \left[ \begin{array}{cc}
    0        & I_l     \\
    - I_l    & 0
         \end{array} \right]
  \end{eqnarray*}
Its Lie algebra ${\mathfrak {sp}} (2l) $ consists of $2l  \times
2l$ matrices of the following form:
  \begin{eqnarray}
  \left[ \begin{array}{ccc}
    a        & b   \\
    c        & -{}^t a
         \end{array} \right]
   \label{eq_matrix}
  \end{eqnarray}
where $a, b, c$ are $l \times l$ matrices, $b, c$ are symmetric.
Let ${\mathfrak b} ({\mathfrak {sp}}(2l))$  be the Borel
subalgebra consisting of matrices of the form (\ref{eq_matrix})
with $c=0$ and $a$ upper triangular. The Cartan subalgebra
${\mathfrak h} ({\mathfrak {sp}} (2l) )$ consists of diagonal
matrices $\text{diag} (t_1, \ldots, t_l, -t_1, \ldots, -t_l )$.
An irreducible representation of $Sp (2l)$ can be parameterized by
its highest weight $\la \in  \Sigma (C) :=\PA^l$.
%
%
%
\subsection{Additional notations}
The notations introduced in this preliminary section are close to
but do not always coincide with \cite{W1}.  For example, our $\binf$
is $\tilde{b}_\infty$ there.

We shall denote by $L(x_\infty; \hL)$, for $x =b, c, d$, the
irreducible $x_\infty$-module of highest weight $\hL$. The level can
be read off from $\Lambda$.

Given a classical Lie group $G$ of type $X$, we shall denote by
$V_\la (G)$ the irreducible $G$-module parameterized by $\la \in
\Sigma(X)$. Let $W(X)$ be the Weyl group of type $X$.

We will denote the roots of the Lie algebra of $G$ by standard
notations $\ep_i \pm \ep_j, \ep_i, 2\ep_i$ etc, and by $(\;\;)$
the bilinear form such that $(\ep_i,\ep_j) =\delta_{ij}$. Let
$\rho$ denote half the sum of positive roots.

\section{Correlation functions on $d_\infty$-modules of level $l$}
\label{sec:dD}

\subsection{The Fock space $\fock^l$}

Let
 $\underline{\Z}$ denote $\hf + \Z \mbox{ or } \Z$,
and set
 $$
  \epsilon =
  \left\{ \begin{array}{ll}
  0, & \text{ if }  \underline{\Z} = \hf+\Z \\
  \hf & \text{ if }  \underline{\Z} = \Z.
  \end{array}
  \right.
 $$
%
Consider a pair of fermionic fields
$$
 \psi^{+}(z) = \sum_{n \in \underline{\Z} } \psi^{+}_n z^{ -n -\hf + \epsilon},
  \quad \psi^{-}(z) = \sum_{n \in \underline{\Z} }
                           \psi^{-}_n z^{ -n -\hf + \epsilon},
 $$
with the following anti-commutation relations
\begin{eqnarray*}
 {[}\psi^{+} _m, \psi^{-}_n ]_{+}
= \delta_{m+n,0},  \qquad
 {[}\psi^{\pm} _m, \psi^{\pm}_n ]_{+} = {0}.
\end{eqnarray*}
Denote by $\mathcal F$ the Fock space of the fermionic fields
$\psi^{\pm}(z)$ generated by a vacuum vector $\vac$ which
satisfies
\begin{eqnarray*}
  \psi^{-}_n \vac = \psi^{+}_n \vac = 0
  & \text{ for } n \in \hf +\Z_+, &\mbox{ if }
   \underline{\Z} = \hf +\Z;                       \\
 \psi^{-}_{n+1} \vac = \psi^{+}_n \vac = 0
 & \text{ for } n \in \Z_{+},
   & \mbox{ if } \underline{\Z} = \Z.
\end{eqnarray*}
We have the standard charge decomposition (cf. \cite{MJD})
$$ \fock =\bigoplus_{k \in \Z} \fock^{(k)}.
$$
Each $\fock^{(k)}$ becomes an irreducible module over a certain
Heisenberg Lie algebra. The shift operator $\textsf{S}:
\fock^{(k)} \rightarrow \fock^{(k+1)}$ matches the highest weight
vectors and commutes with the creation operators in the Heisenberg
algebra.

Now we take $l$ pairs of fermionic fields, $\psi^{\pm,p} (z) \;( p
= 1, \dots, l)$ and denote the corresponding Fock space by
$\fock^l$. Introduce the following generating series
\begin{eqnarray}
  E(z,w) & \equiv & \sum_{i,j \in \Z} E_{ij} z^{i-1 + 2\epsilon }
     w^{-j  }
   = \sum_{p =1}^l \no\psi^{+,p} (z) \psi^{-,p} (w)\no,  \label{eq_genelin}
\end{eqnarray}
where the normal ordering $::$ means that the operators
annihilating $\vac$ are moved to the right with a sign. It is well
known that the operators $E_{ij} \; (i,j \in \Z)$ generate a
representation in $\fock^l$ of the Lie algebra $\hgl$ with level
$l$.

Let
\begin{eqnarray} \label{e-}
e^-_{pq}  =\sum_{r\in
\underline{\Z}}\no\psi_{r}^{-p}\psi^{-q}_{-r}\no,\quad
e_{pq}^+ = \sum_{r\in\underline{\Z}}
\no\psi_{r}^{+p}\psi^{+q}_{-r}\no, \quad p \neq q,
\end{eqnarray}
 and let
\begin{eqnarray}  \label{e+}
e_{pq}  =\sum_{r\in\underline{\Z}}
\no\psi_{r}^{+p}\psi^{-q}_{-r}\no +\delta_{pq} \eps.
\end{eqnarray}
The operators $e^+_{pq}, e_{pq}, e^-_{pq} \;(p, q = 1, \cdots, l)$
generate Lie algebra $\mathfrak {so}(2l)$ (cf. \cite{FF, W1}).

\subsection{The  $(O(2l), \dinf)$-Howe duality} \label{sec:HoweDd}

Now let $\underline{\Z} = \hf+\Z$ for the remainder of
Section~\ref{sec:dD}.

The representation of the Lie algebra $\dinf$ on $\fock^l$ is
given by (cf. \cite{DJKM2})
\begin{eqnarray} \label{gs:donf}
&& \sum_{i,j\in\Z}(E_{i,j} -E_{1-j,1-i})z^{i-1}w^{-j} \nonumber \\
 &=& \sum_{p =1}^l
\no\psi^{+,p}(z) \psi^{-,p}(w)\no -\no\psi^{+,p}(w)\psi^{-,p}(z)
\no
\end{eqnarray}

The action of $\mathfrak {so}(2l)$ can be integrated to the action
of the Lie group $SO(2l)$ on $\fock^l$. In particular, the
operators $e_{pq}\; ( p, q = 1, \cdots, l)$ form a Lie subalgebra
${\mathfrak {gl}}(l)$. We identify the Borel subalgebra
${\mathfrak b}( {\mathfrak {so}}(2l) )$ with the one generated by
$e^+_{pq} \;(p \neq q), e_{pq}\; (p \leq q), p, q = 1, \ldots, l.$
Note that $\tau \in O(2l)$ defined in (\ref{eq_tau}) commutes with
the action of $\dinf$ on $\fock^l$. The following lemma summarizes
Lemmas 3.2, 3.3 in \cite{W1}.
\begin{lemma}
The action of the Lie group $O(2l)$ commutes with the action of
$\dinf$ on $\fock^l$.
\end{lemma}

We define a map $\Lambda: \Sigma (D) \longrightarrow {\dinf}_0^*$
by sending $ \lambda = (m_1, \cdots, { m}_l)$, where $m_l
> 0,$ to
$$  \Lambda (\lambda) =
   (l -i ) \hL_0^d + ( l -i ) \hL_1^d + \sum_{k =1}^i \hL_{m_k}^d, $$
sending $(m_1, \cdots, m_j, 0, \ldots, 0)$, where $j <l,$ to
$$  \Lambda (\lambda) =
   (2l -i -j) \hL_0^d + (j-i ) \hL_1^d + \sum_{k =1}^i \hL_{m_k}^d, $$
and sending $ (m_1, \ldots, m_j, 0, \ldots, 0 ) \otimes {\det}$,
where $j <l,$  to
 $$  \Lambda (\lambda) =
  (j-i) \hL_0^d + (2l -i -j) \hL_1^d + \sum_{k =1}^i \hL_{m_k}^d, $$
if $m_1 \geq \ldots m_{i} > m_{i+1} = \ldots = m_{j} =1
   > m_{j+1} = \ldots = m_l = 0.$

\begin{prop} \cite[Theorem 3.2]{W1} \label{prop:HoweDd}
We have the following decomposition of $(O(2l),d_\infty)$-modules:
\begin{eqnarray}\label{duality}
\fock^{l} \cong \bigoplus_{\la\in\Sigma(D)} V_\la(O(2l))\otimes
L(d_\infty;\Lambda(\la)).
\end{eqnarray}
\end{prop}
\subsection{The main results of \cite{BO, CW}}

Let $t$ be an indeterminant. We define the following operators in
$d_\infty$ (cf. \cite{W1}):
\begin{eqnarray*}
 \no\D(t)\no &=& \sum_{k\in \mathbb N}(t^{k-\hf} -
t^{\hf-k})(E_{k,k}-E_{1-k,1-k}), \\
 \D(t) &=& \no \D(t)\no +\frac{2}{t^\hf-t^{-\hf}}C.
\end{eqnarray*}
When acting on $\fock^l$, these operators can be written in terms
of the operators $\psi^\pm_n$ by (\ref{gs:donf}) as
\begin{eqnarray*}
\no \D(t)\no &=&\sum_{p =1}^l \sum_{k\in\HZ} t^k (\no
\psi^{+,p}_{-k} \psi^{-,p}_{k}\no + \no\psi^{-,p}_{-k}
\psi^{+,p}_k \no), \\
 \D(t) &=& \sum_{p =1}^l \sum_{k\in\HZ} t^k
(\psi^{+,p}_{-k} \psi^{-,p}_{k} +\psi^{-,p}_{-k} \psi^{+,p}_k).
\end{eqnarray*}

Recall that Bloch and Okounkov \cite{BO} introduced  the following
operators in $\hgl$
$$\no\A(t)\no
= \sum_{k\in\Z} t^{k-\hf}E_{k,k}, \quad
\A(t) = \no\A(t)\no + \frac{1}{t^\hf-t^{-\hf}} C.
$$
We easily verify that
\begin{equation} \label{eq:AD}
\D(t) = \A(t) - \A(t^{-1}), \qquad \no \D(t)\no = \no
\A(t)\no - \no \A(t^{-1})\no.
\end{equation}

Given $\la =(m_1,\ldots,m_l)\in \Sigma(A)$, we denote by $\hL(\la)$
the $\hgl$-highest weight $\hL^a_{m_1}+\cdots+\hL^a_{m_l}$. The
energy operator $L_0$ on the $\hgl$-module $L(\hgl;\hL(\la))$ with
highest weight vector $v_{\hL(\la)}$ is characterized by
\begin{eqnarray}
L_0 \cdot v_{\hL (\la)}
 &=& \hf \Vert\la\Vert^2  \cdot v_{\hL (\la)}, \label{weight} \\
 {[}L_0, E_{ij}] &=& (i-j) E_{ij}, \nonumber
\end{eqnarray}
where $$\Vert\la\Vert^2 :=\la_1^2 +\la_2^2 +\cdots + \la_l^2,$$
On $\fock^l$, we can realize $L_0$ as
$$L_0 =
\sum_{p=1}^l\sum_{k\in\Z+\hf}
k\no\psi_{-k}^{+,p}\psi_{k}^{-,p}\no.$$

The {\em $n$-point $\hgl$-correlation function of level $l$}
associated to $\la$ is defined in \cite{BO} for $l=1$ and in
\cite{CW} for general $l$ as
\begin{equation*}
\mathfrak A_\la^l (q;{\bf t}) \equiv
 \mathfrak A_\la^l (q;t_1,\ldots,t_n)
 :=\trace_{L(\hgl,\hL(\la))}(q^{L_0} \A(t_1)\A(t_2)\cdots
 \A(t_n)).
\end{equation*}
Here and below we denote ${\bf t} =(t_1, \ldots, t_n)$.

Let $ (a; q)_\infty := \prod_{r=0}^\infty (1-aq^r)$. Define the
theta function
\begin{eqnarray}
 \label{theta:def} \Theta (t) &:=&
 (t^{\hf} -t^{-\hf})(q;q)_\infty^{-2} (qt; q)_\infty(qt^{-1};q)_\infty \\
 \label{theta:deriv} \Theta^{(k)} (t) &:=&
 \left(t\frac{d}{dt} \right)^k \Theta (t), \quad\text{ for }k\in\Z_+.
\end{eqnarray}
Denote by $F_{bo}(q; {\bf t})$ or $F_{bo}(q;t_1,\ldots,t_n)$ the
following expression
\begin{eqnarray}  \label{eq:bo}
\frac{1}{(q; q)_\infty}\cdot \sum_{\sigma\in S_n} \frac{{\rm det}
\Big{(}\frac{\Theta^{(j-i+1)}(t_{\sigma(1)}\cdots
t_{\sigma(n-j)})}{(j-i+1)!} \Big{)}_{i,j=1}^n}
{\Theta(t_{\sigma(1)}) \Theta(t_{\sigma(1)}t_{\sigma(2)})\cdots
\Theta(t_{\sigma(1)}t_{\sigma(2)}\cdots t_{\sigma(n)})}.
\end{eqnarray}
It is understood here that $1/ (-k)!=0$ for $k>0$, and for $n=1$,
we have $F_{bo}(q; t) = (q;q)_\infty^{-1} \Theta(t)^{-1}$. The
following summarizes the main results of Bloch-Okounkov \cite{BO}
for $l=1$ and Cheng-Wang \cite{CW} for general $l \ge 1$.
\begin{theorem}  \label{th:CWmain}
Associated to $\la =(\la_1,\ldots,\la_l)$, where $\la_1\ge\ldots
\ge \la_l$ and $\la_i \in\Z$, the $n$-point $\hgl$-function of
level $l$ is given by
\begin{eqnarray*}
\mathfrak A^l_\la(q; {\bf t}) =
 q^{\frac{\Vert\la\Vert^2}{2}} (t_1t_2\cdots t_n)^{|\la|} \prod_{1\le i<j\le
l} (1-q^{\la_i-\la_j+j-i}) \cdot
 F_{bo}(q;{\bf t})^l
\end{eqnarray*}
where $|\la| := \la_1 +\cdots +\la_l.$
\end{theorem}
In the simplest case, i.e. $l =n =1$, we have
$$
\mathfrak A^1_\la(q; t) =q^{\frac{\la^2}{2}}t^\la\cdot F_{bo}(q; t)
=\frac{q^{\frac{\la^2}{2}}t^\la}{(q;q)_\infty \Theta(t)}.
$$

\subsection{The $1$-point $\dinf$-functions of level $l$}

\begin{definition}
The $n$-point $\dinf$-correlation function of level $l$ associated
to $\la =(\la_1, \ldots, \la_l) \in \PA^l$ is
\begin{eqnarray*}
&& \mathfrak D^l_\la(q,{\bf t}\,) \equiv \mathfrak
D^l_\la(q,t_1,\ldots, t_n) \\
&:=& \left \{
\begin{array}{ll}
  \trace_{L(d_\infty;\Lambda(\la))} q^{L_0}\D(t_1)\cdots \D(t_n),
\qquad\;\; \;\quad \qquad \text{ if } \la_l \neq 0,
\\
\trace_{L(d_\infty;\Lambda(\la)) \oplus
L(d_\infty;\Lambda({\la}\otimes \det))} q^{L_0}\D(t_1)\cdots
\D(t_n), \text{ if } \la_l = 0.
\end{array}
\right.
\end{eqnarray*}
(The operator $L_0$ is defined in the same way as for a
$\hgl$-module.)
\end{definition}

\begin{remark} \label{rem:Oinf}
A justification of this definition when $\la_l = 0$ is as follows.
The two weights $\Lambda(\la)$ and $\Lambda(\la \otimes \det)$ are
interchanged by a $\dinf$-Dynkin diagram automorphism. Thus, the
direct sum $L(d_\infty;\Lambda(\la)) \oplus L(d_\infty;\Lambda(\la
\otimes \det))$ can be regarded as an irreducible module of the
orthogonal group associated to $\dinf$.
\end{remark}

In this subsection, we will restrict to $n=1$ for notational
simplicity.

\begin{theorem} \label{th:1ptD}
The $1$-point $\dinf$-function of level $l$ is given by
$$
\mathfrak D^l_\la(q,t) = {\frac1{(q;q)_\infty^l \Theta(t)^l}}
\cdot
\sum_{\sigma\in W(D_l)}(-1)^{\ell(\sigma)} q^{\frac{\Vert
\la+\rho-\sigma(\rho) \Vert^2}{2}}
\prod_{a=1}^l(t^{k_a}+t^{-k_a})$$
where $k_a = (\la+\rho-\sigma(\rho), \ep_a)$.
\end{theorem}

We first prepare a few lemmas for the proof of this theorem. By a
character of a module $V$ of $SO(2l)$ or $O(2l)$, we mean
$\trace_{V} (z_1^{e_{11}}\ldots z_l^{e_{ll}})$. Let $|a_{ij}|$
denote the determinant of a matrix $(a_{ij})$.

\begin{lemma} \label{lem:charD}
Denote by $\text{ch}^o_\la (z_1,\dots,z_l)$ the character of the
irreducible $O(2l)$-module $V_\la(O(2l))$. If $\la_l = 0$, then
\begin{eqnarray*}
 \text{ch}^o_\la (z_1,\dots,z_l) =
\frac{\left |z_j^{\la_i+l-i} + z_j^{-(\la_i+l-i)}\right
|}{\left|z_j^{l-i}+z_j^{-(l-i)}\right |}.
\end{eqnarray*}
If $\la_l \neq 0$, then
$$\text{ch}^o_\la (z_1,\dots,z_l) = \frac{2\left|z_j^{\la_i+l-i} +
z_j^{-(\la_i+l-i)}\right|}{\left|z_j^{l-i}+z_j^{-(l-i)}\right|}.$$
\end{lemma}

\begin{proof}
The character of an irreducible $SO(2l)$-module $V_\la(SO(2l))$,
denoted by $\text{ch}^{so}_\la$, is well known to be as follows
(cf. \cite[pp. 410]{FH}):
\begin{eqnarray} \label{char:so}
\text{ch}^{so}_\la (z_1,\ldots, z_l)
 =
\frac{\left|z_j^{\la_i+l-i} + z_j^{-(\la_i+l-i)}\right| +
\left|z_j^{\la_i+l-i} -
z_j^{-(\la_i+l-i)}\right|}{\left|z_j^{l-i}+z_j^{-(l-i)}\right|}.
\end{eqnarray}
Recall that an irreducible $O(2l)$-module can be a sum of two
irreducible $SO(2l)$-modules or remain to be irreducible as a
$SO(2l)$-module, depending on whether $\la_l$ is nonzero or not.
If $\la_l = 0$, the second determinant in the numerator of
(\ref{char:so}) vanishes and hence the character formula for
$\text{ch}^o_\la$ of $V_\la(O(2l))$ follows from
$\text{ch}^{so}_\la$. If $\la_l \neq 0$, then $V_\la(O(2l))=
V_\la(SO(2l)) \oplus V_{\overline{\la}} (SO(2l))$. Note that the
second determinant terms in the numerators of $\text{ch}^{so}_\la$
and $\text{ch}^{so}_{\overline{\la}}$ (cf. (\ref{char:so})) are
opposite to each other. Now the formula for $\text{ch}^o_\la$
follows.
\end{proof}

\begin{lemma} \label{lem:sym}
We have
\begin{eqnarray*}
\trace_{\fock^{(k)}}q^{L_0}\D(t) &=&
\frac{1}{2}(t^k+t^{-k})q^{\frac{k^2}{2}} \trace_{\fock^{(0)}} q^{L_0}\D(t), \quad k\in \Z,\\
\trace_\fock  z^{e_{11}}q^{L_0}\D(t)
 &=& \trace_{\fock^{(0)}}
q^{L_0}\D(t) \cdot
\sum_{k\in\Z}\frac{(t^k+t^{-k})}{2}z^kq^{\frac{k^2}{2}}.
\end{eqnarray*}
\end{lemma}

\begin{proof}
These two identities are clearly equivalent, and it suffices to
prove the first one. We refer to \cite{MJD} or \cite[Appendix
$A$]{Ok} for properties of the shift operator $\textsf{S}$ on
$\fock$. Then,
$$\begin{aligned} \trace_{\fock^{(k)}} q^{L_0}\D(t)&=
\trace_{\fock^{(0)}}\textsf{S}^{-k}q^{L_0}(\A(t)-\A(t^{-1}))\textsf{S}^k
\\ &=
\trace_{\fock^{(0)}}q^{L_0+\frac{k^2}{2}}(t^k\A(t)-t^{-k}\A(t^{-1}))
\\ &= q^{\frac{k^2}{2}}\trace_{\fock^{(0)}}
q^{L_0}(t^k\A(t)-t^{-k}\A(t^{-1}))
\end{aligned}
$$
and
$$\begin{aligned}
\trace_{\fock^{(-k)}} q^{L_0}\D(t) &=
\trace_{\fock^{(0)}}\textsf{S}^{k}q^{L_0}(\A(t)-\A(t^{-1}))\textsf{S}^{-k}
\\ &=
\trace_{\fock^{(0)}}q^{L_0+\frac{k^2}{2}}(t^{-k}\A(t)-t^{k}\A(t^{-1}))
\\ &= q^{\frac{k^2}{2}}\trace_{\fock^{(0)}}
q^{L_0}(t^{-k}\A(t)-t^{k}\A(t^{-1})).
\end{aligned}
$$
By (\ref{eq:AD}), $t^k\A(t)-t^{-k}\A(t^{-1}) +t^{-k}\A(t)
-t^k\A(t^{-1}) = (t^k +t^{-k}) \D(t)$.  It follows that
$$\trace_{\fock^{(-k)}}q^{L_0}\D(t) +
\trace_{\fock^{(k)}}q^{L_0}\D(t) =(t^k+t^{-k})q^{\frac{k^2}{2}}
\trace_{\fock^{(0)}} q^{L_0}\D(t).
$$

Since $\fock^{(-k)}$ and $\fock^{(k)}$ are isomorphic as
$\dinf$-modules and $\D(t)$ lies in $\dinf$, we have
$\trace_{\fock^{(-k)}}q^{L_0}\D(t) =
\trace_{\fock^{(k)}}q^{L_0}\D(t)$, and the result
follows.\end{proof}

\begin{lemma} \label{lem:AA=D}
 We have $\trace_{\fock^{(0)}} q^{L_0}\D(t) = 2 F_{bo}(q,t) =\frac2{(q;q)_\infty
\Theta(t)}$.
\end{lemma}

\begin{proof}
As a special case of Theorem~\ref{th:CWmain}, we note that
$$\mathfrak A^1_{(0)}(q; t) =\trace_{\fock^{(0)}}
q^{L_0}\A(t) =F_{bo}(q;{t}) =\frac1{(q,q)_\infty \Theta(t)}.
$$
In particular, $\trace_{\fock^{(0)}} q^{L_0}\A(t^{-1})
=-\trace_{\fock^{(0)}} q^{L_0}\A(t)$. Thus, it follows by
(\ref{eq:AD}) that $\trace_{\fock^{(0)}} q^{L_0}\D(t)
=2F_{bo}(q;{t}).$
\end{proof}

Set
$$
c_\la =\left\{
\begin{array}{ll}
 1, & \text{ if }\la_l =0 \\
 2, & \text{ if }\la_l \neq 0.
\end{array}
\right.
$$

We have the following lemma similar to \cite[Lemma 6]{CW}.

\begin{lemma} \label{lem:qseries}
We have the following identity:
$$\prod_{i=1}^l \trace_\fock z_i^{e_{ii}}q^{L_0}\D(t)
 =
 \sum_\la c_\la \cdot \frac{\left|z_j^{\la_i+l-i} +
z_j^{-(\la_i+l-i)}\right|}{\left |z_j^{l-i}+z_j^{-(l-i)}\right|}
\cdot \mathfrak D^l_\la(q,t).$$
\end{lemma}

\begin{proof}
It follows by computing the trace of $z_1^{e_{11}}\cdots
z_l^{e_{ll}} q^{L_0}\D(t)$ of both sides of the
$(O(2l),\dinf)$-duality (\ref{duality}). Note that the factor $2$ in
$\text{ch}^o_\la$ in Lemma~\ref{lem:charD} for $\la_l \neq 0$ gives
rise to the factor $c_\la$.
\end{proof}

Let ${\bf z}^\mu$ denote $z_1^{\mu_1}\cdots z_l^{\mu_l}$ for
$\mu=(\mu_1,\ldots,\mu_l)$.
\begin{lemma} \label{weylcharindep}
Among all the monomials ${\bf z}^\mu$ in the expansion of the
determinant ${\left|z_j^{\la_i+l-i} + z_j^{-(\la_i+l-i)}\right|}$,
there is exactly one ``dominant" monomial with $\mu_1 \ge \ldots
\ge \mu_l \ge 0$, that is, ${\bf z}^{\la+\rho}=\prod_{i=1}^l
z_i^{\la_i+l-i}$. Its coefficient is equal to $2/ c_\la$, or more
concretely, it is $1$ if $\la_l \neq 0$ and $2$ if $\la_l =0$.
\end{lemma}

\begin{proof}
The first part is clear by inspection. Note that the coefficient
$2$ when $\la_l =0$ comes from the last row of the determinant.
\end{proof}

\begin{proof}[Proof of Theorem~\ref{th:1ptD}]

It follows by Lemmas~\ref{lem:sym}, \ref{lem:AA=D}, and
\ref{lem:qseries} that
\begin{eqnarray} \label{eq:iden}
&& F_{bo}(q,t)^l  \cdot
  \prod_{a=1}^l \left(\sum_{k_a\in\Z}
 {(t^{k_a} +t^{-{k_a}})} z_a^{k_a} q^{\frac{k_a^2}{2}} \right)
\cdot
 \hf {\left |z_j^{l-i}+z_j^{-(l-i)}\right|} \nonumber \\
 &=& \sum_\la  {c_\la}/2 \cdot {\left|z_j^{\la_i+l-i} +
z_j^{-(\la_i+l-i)}\right|} \cdot \mathfrak D^l_\la(q,t).
\end{eqnarray}

Recall (\cite[(24.38)]{FH}) that the Weyl denominator of type
$D_l$ is
\begin{eqnarray} \label{eq:denom}
\hf {\left |z_j^{l-i}+z_j^{-(l-i)}\right|} = \sum_{\sigma\in
W(D_l)} (-1)^{\ell(\sigma)} {\bf z}^{\sigma(\rho)}.
\end{eqnarray}
The theorem follows when we apply Lemma~\ref{weylcharindep} to
compare the coefficients of the monomial $\prod_{i=1}^l
z_i^{\la_i+l-i}$ on both sides of (\ref{eq:iden}) with
(\ref{eq:denom}) plugged in.
\end{proof}

\subsection{The $n$-point $\dinf$-functions of level $l$}

We now compute the $n$-point $\dinf$-function of level $l$. The
$1$-point calculation in the previous subsection carries over for
general $n$ after suitable modification. Let
\begin{eqnarray} \label{eq:traceF}
\bff(z,q; t_1,\dots,t_n) = \trace_\fock
z^{e_{11}}q^{L_0}\D(t_1)\cdots \D(t_n)
\end{eqnarray}

The following lemma is the $n$-point generalization of
Lemma~\ref{lem:sym}.

\begin{lemma} \label{lemma:nptref}
We have
\begin{eqnarray*}
 \bff(z,q; t_1,\dots,t_n) = \sum_{k\in\Z} z^k q^{\frac{k^2}{2}}
 \sum_{\vec{\eps}  \in \{\pm 1\}^n}  [\vec{\eps}] \cdot
 \left(\Pi{\bf t}^{\vec{\eps}} \right)^k \cdot
 F_{bo}(q;{\bf t}^{\vec{\eps}})
\end{eqnarray*}
where we denote $\vec{\eps} =(\eps_1,\eps_2,\ldots,\eps_n)$,
$[\vec{\eps}] =\eps_1\eps_2 \cdots \eps_n$, $\Pi{\bf
t}^{\vec{\eps}} =t_1^{\eps_1}\cdots t_n^{\eps_n}$,
 and $ F_{bo}(q;{\bf
t}^{\vec{\eps}}) = F_{bo}(q;t_1^{\eps_1},\dots,t_n^{\eps_n})$.
\end{lemma}

\begin{proof}
We calculate using (\ref{eq:AD}) that
$$
\begin{aligned}
\D(t_1)\cdots \D(t_n)
 &= \prod_{j=1}^n \left(\A(t_j) -\A(t_j^{-1})\right) \\
& = \sum_{\vec{\eps}\in\{\pm 1\}^n} \eps_1\eps_2 \cdots \eps_n
\A(t_1^{\eps_1}) \A(t_2^{\eps_2})\cdots \A(t_n^{\eps_n})
\end{aligned}
$$
It is known (cf. \cite{Ok, CW}) by the same type of argument as in
Lemma~\ref{lem:sym} that
$$\trace_{\fock} z^{e_{11}}q^{L_0}\A(t_1)\cdots \A(t_n) =
\sum_{k\in\Z} z^k q^{\frac{k^2}{2}} (t_1\cdots t_n)^k
F_{bo}(q,t_1, \ldots, t_n).
$$
Now the lemma follows.
\end{proof}

\begin{lemma} \label{lem:q2}
We have the following identity:
$$\prod_{i=1}^l \bff(z_i,q;t_1,\ldots,t_n)
 =  \sum_\la
 c_\la \cdot \frac{\left|z_j^{\la_i+l-i} +
z_j^{-(\la_i+l-i)}\right|}{\left |z_j^{l-i}+z_j^{-(l-i)}\right|}
\cdot \mathfrak D^l_\la(q,t_1,\ldots,t_n).$$
\end{lemma}

\begin{proof}
The lemma is a straightforward $n$-point version of
Lemma~\ref{lem:qseries}, which is proved in the same way as
before.
\end{proof}

\begin{theorem} \label{th:nlD}
The $n$-point $\dinf$-correlation function of level $l$ associated
to $\la \in \PA^l$ is given by
\begin{eqnarray*}
&& \mathfrak D_\la^l(q;t_1,\dots,t_n) =\\
&& \sum_{\sigma\in W(D_l)}\left(-1\right)^{\ell(\sigma)}
q^{\frac{\Vert \la+\rho-\sigma(\rho) \Vert^2}{2}}
 \prod_{a=1}^l
 \Big(
\sum_{\vec{\eps}_a\in\{\pm 1\}^n}
 [\vec{\eps}_a] (\Pi {\bf t}^{\vec{\eps}_a})^{k_a}
 F_{bo}(q;{\bf t}^{{\vec{\eps}_a}})
 \Big)
\end{eqnarray*}
where $k_a = (\la+\rho-\sigma(\rho), \ep_a)$.
\end{theorem}

\begin{proof}
The proof is the same as the proof of Theorem~\ref{th:1ptD} by
using now Lemmas~\ref{lemma:nptref} and \ref{lem:q2}.
\end{proof}
\subsection{A refined $1$-point function of level $1$}

Denote by $\fock^{(0)}_\pm$ the $(\pm 1)$-eigenspace of $\tau \in
O(2)$ acting on $\fock^{(0)}$ where we recall that
$$
\tau = \begin{bmatrix} 0 & 1 \\ 1 & 0 \end{bmatrix}.
$$
As an $(O(2),d_\infty)$-module, we have
$$
\fock = \oplus_{m\in\mathbb N} (\fock^{(m)}\oplus
\fock^{(-m)})\oplus \fock^{(0)}_+\oplus\fock^{(0)}_-.
$$
As $\dinf$-modules, all $\fock^{(m)}, \fock^{(-m)}, \fock^{(0)}_+,
\fock^{(0)}_-$ are irreducible, and moreover, $\fock^{(m)} \cong
\fock^{(-m)}$.
 The main result of this subsection is the following.

\begin{theorem} \label{th:refined}
The $\trace_{\fock^{(0)}_\pm}q^{L_0}\D(t)$, which we refer to as a
refined $1$-point $\dinf$-correlation function, is equal to
\begin{align*}
& \frac{1}{(q;q)_\infty\Theta(t)} \\
&  \quad \pm
 (q;q^2)_\infty
 \left ( \sum_{r=0}^\infty
\Big(\frac{q^{r+1}t^{-\hf}}{1-q^{2(r+1)}t^{-1}} -
\frac{q^{r+1}t^\hf} {1-q^{2(r+1)}t}\Big) +
\frac{1}{t^\hf-t^{-\hf}} \right)
\end{align*}
which is equivalent to
\begin{align*}
& \frac{1}{(q;q)_\infty\Theta(t)}
 \pm
 (q;q^2)_\infty
 t\frac{d}{dt}\ln\left ( \frac{(-t^{-\hf};q)_\infty
 (-qt^{\hf};q)_\infty}{(t^{-\hf};q)_\infty (qt^{\hf};q)_\infty} \right).
\end{align*}
\end{theorem}

We denote
\begin{eqnarray*}
G(t) &=& \trace_{\fock^{(0)}_+}q^{L_0}\D(t)-\trace_{\fock^{(0)}_-}
q^{L_0}\D(t)
\\
\no G(t)\no &=& \trace_{\fock^{(0)}_+}q^{L_0}\no\D(t)\no
-\trace_{\fock^{(0)}_-} q^{L_0} \no\D(t)\no.
\end{eqnarray*}
To compute $\trace_{\fock^{(0)}_\pm}q^{L_0}\D(t)$ it suffices to
compute their difference $G(t)$, since $\mathfrak D^1_{(0)} (t)
=\trace_{\fock^{(0)}_+}q^{L_0}\D(t) +\trace_{\fock^{(0)}_-}
q^{L_0}\D(t)$ has been calculated.

Recall that the {\em rank} $\text{rk} (\la)$ of a partition $\la
=(\la_1,\la_2,\ldots)$ is the cardinality of the set $\{i \mid
\la_i \geq i\}$.

\begin{lemma} \label{lem:G}
We have
$$
\no G(t)\no  = 2 \sum_{\la = \la^t} (-1)^{\text{rk}(\la)}q^{|\la|}
\sum_i\left(t^{\la_i-i+1/2} -t^{-(\la_i-i+1/2)}\right).
$$
\end{lemma}

\begin{proof}
Note that $\tau$ sends $\psi^+_n$ to $\psi^-_n$ for each $n$ and
vice-versa. Thus, $ \fock^{(0)}_+$ has a basis given by
\begin{eqnarray*}
&& \{\psi^-_{m_1}\cdots\psi^-_{m_r}\psi^+_{n_1}
\cdots\psi^+_{n_r}|0\rangle
 + \psi^-_{n_1}\cdots\psi^-_{n_r}\psi^+_{m_1}\cdots\psi^+_{m_r}|0\rangle \\
&& \quad \mid m_1 < \cdots < m_r < 0, n_1 < \cdots < n_r < 0, m_i
\not=
n_i\, \text{ for some } i\} \\
&  \cup&
\{\psi^-_{m_1}\cdots\psi^-_{m_r}\psi^+_{m_1}\cdots\psi^+_{m_r}|0\rangle\st
\,m_1 < \cdots < m_r < 0,  r  \text{ even}\}
\end{eqnarray*}
and $\fock^{(0)}_-$ has a basis given by
\begin{eqnarray*}
&& \{\psi^-_{m_1}\cdots\psi^-_{m_r}\psi^+_{n_1} \cdots\psi^+_{n_r}
|0\rangle -\psi^-_{n_1}\cdots\psi^-_{n_r}\psi^+_{m_1}
\cdots\psi^+_{m_r} |0\rangle \\
&& \quad \mid m_1 < \cdots < m_r < 0, n_1 < \cdots < n_r < 0, m_i
\not= n_i \text{ for some }i \} \\
&\cup&  \{\psi^-_{m_1}\cdots\psi^-_{m_r}\psi^+_{m_1}
\cdots\psi^+_{m_r} |0\rangle \st \,m_1 < \cdots < m_r < 0, r\,{\rm
odd}\}
\end{eqnarray*}
$\fock^{(0)}_+$ (resp. $\fock^{(0)}_-$) has highest weight vector
$|0\rangle$ (resp. $\psi^-_{-\hf}\psi^+_{-\hf}|0\rangle$).

The action of $\no \D(t)\no$ on $\fock^{(0)}$ can be described
explicitly:
$$\begin{aligned}
\no \D(t)\no\, &\psi^-_{m_1}\cdots\psi^-_{m_r}\psi^+_{n_1}\cdots
\psi^+_{n_r}|0\rangle
\\ &=  \sum_i \left( t^{-n_i}-t^{n_i}
 + t^{-m_i}-t^{m_i}\right)
 \psi^-_{m_1}\cdots\psi^-_{m_r}\psi^+_{n_1}
\cdots\psi^+_{n_r} |0\rangle.
\end{aligned}
$$

It is well known that $\fock^{(0)}$ can be identified with an
irreducible module of Heisenberg algebra and its basis is
parameterized by partitions. Given an element
$\psi^-_{-p_1}\cdots\psi^-_{-p_r}\psi^+_{-q_1}\cdots
\psi^+_{-q_r}|0\rangle \in \fock^{(0)}$, where $p_1>\ldots>p_r>0,
q_1>\ldots>q_r>0$, the indices $(p_1,\ldots, p_r \mid
q_1,\ldots,q_r)$ are exactly the Frobenius coordinates of a
partition $\la$ (which by our convention here uses half-integers).
It is well known that (cf. e.g. \cite[Lemma 5.1]{BO})
$$\sum_{i=1}^l \left(t^{\la_i-1+\hf}-t^{-i+\hf}\right) =
\sum_{k=1}^r\left(t^{p_k}-t^{-(-q_k)}\right).$$

We then compute that
$$
\begin{aligned}
\trace_{\fock^{(0)}_+} q^{L_0} \no\D(t)\no
 = & X + 2\mathop{\sum_{\la = \la^t}}_{\text{rk}(\la) \text{ even}}
q^{|\la|}\sum_i\left(t^{\la_i-i+1/2}-t^{-(\la_i-i+1/2)}\right)
\\
\trace_{\fock^{(0)}_-} q^{L_0} \no\D(t)\no =
 &X +2\mathop{\sum_{\la = \la^t}}_{\text{rk}(\la) \text{ odd}}
q^{|\la|}\sum_i\left(t^{\la_i-i+1/2} -t^{-(\la_i-i+1/2)} \right),
\end{aligned}
$$
for the same $X$ whose explicit form is irrelevant here. Now $\no
G(t)\no$ is given by the difference of the above two formulas.
\end{proof}

\begin{prop}  \label{prop:Gt}
We have
$$
\begin{aligned}\no G(t)\no &=
2(q;q^2)_\infty\cdot\sum_{n=1}^\infty
\frac{q^{2n-1}\left(t^{-n+\hf}-t^{n-\hf}\right)}{1-q^{2n-1}}\\
&= 2(q;q^2)_\infty\cdot\sum_{r=0}^\infty
\left(\frac{q^{r+1}t^{-\hf}}{1-q^{2(r+1)}t^{-1}} -
\frac{q^{r+1}t^\hf} {1-q^{2(r+1)}t}\right).
\end{aligned}$$
and
$$
G(t) = \no G(t)\no + (q;q^2)_\infty \frac{2}{t^\hf-t^{-\hf}}.
$$
\end{prop}

\begin{proof}
There is a canonical bijection between the set of symmetric
partitions (i.e. $\la$ such that $\la^t =\la$) and the set $\OSP$ of
odd strict partitions. The bijection is achieved by setting the
parts of a new partition $\mu$ to be the hook lengths of the
diagonal nodes from the original symmetric partition $\la$ (i.e.
$\mu_i =2\la_i-2i+1$). Here is an example:
$$
(5,4,3,2,1) \mapsto (9,5,1).
$$
Under such a bijection which sends a symmetric partition $\la$ to
$\mu \in \OSP$, we have $rk(\la) = \ell(\mu), |\la| = |\mu|,$ and
$\la_i-i+\hf = \mu_i/2$. By Lemma~\ref{lem:G}, we have a
reformulation
\begin{eqnarray} \label{eq:Gt}
\no G(t)\no
 = 2 \sum_{\mu\in\OSP} (-1)^{\ell(\mu)} q^{|\mu|} \sum_{k\geq 1}
\left(t^{\mu_k/2} - t^{- \mu_k/2}\right).
\end{eqnarray}

By Theorem 8 of \cite{W2} (replacing $q$ therein by $q^2$), we
have the following identity:
$$
\sum_{\mu\in\OSP} \left(z^{\ell(\mu)} q^{|\mu|} \sum_{k\geq 1}
t^{\mu_k/2}\right)
 = (-qz;q^2)_\infty\left(1+\sum_{n=1}^\infty
\frac{q^{2n-1}t^{n-\hf}z}{1+q^{2n-1}z }\right).
$$
Using this identity twice (with the specialization $z=-1$), we
obtain from (\ref{eq:Gt}) the first formula for $\no G(t)\no$. The
second formula for $\no G(t)\no$ follows from the following
identity:
$$
\sum_{n=1}^\infty \frac{zq^{2n-1}t^{n-\hf}}{1+q^{2n-1}z} =
\sum_{r=0}^\infty
\frac{(-1)^rz^{r+1}q^{r+1}t^\hf}{1-q^{2(r+1)}t}.
$$
(This identity follows quickly by expanding the left side as a
power series, interchanging summations, and then summing up.)

Along the same line as the proof of Lemma~\ref{lem:G} and
(\ref{eq:Gt}), we can show that
$$
\trace_{\fock^{(0)}_+} q^{L_0} -\trace_{\fock^{(0)}_-} q^{L_0}
=\sum_{\mu \in \OSP} (-1)^{\ell(\mu)} q^{|\mu|} = (q;q^2)_\infty$$
where the last equation is the specialization at $z=-1$ of the
identity
$$
\sum_{\mu \in \OSP} z^{\ell(\mu)} q^{|\mu|}= (-qz;q^2)_\infty.
$$
Now the formula for $G(t)$ follows from this consideration and
(\ref{eq:AD}).
\end{proof}

\begin{remark}
The function $G(t)$ is essentially a ``super version" of the
$1$-point $\dinf$-correlation function of level $\hf$ (cf.
\cite[Theorem~9]{W2}).
\end{remark}

\begin{proof}[Proof of Theorem~\ref{th:refined}]
We write
$$
\trace_{\fock^{(0)}_\pm}q^{L_0}\D(t)
 = \hf \left( \mathfrak D^1_{(0)} (t)\pm G(t) \right).
 $$

It is known by Theorem~\ref{th:1ptD} that
$$
\mathfrak D^1_{(0)} (t) =\frac{2}{(q;q)_\infty\Theta(t)}.
$$
Now the first formula of the theorem follows from
Proposition~\ref{prop:Gt}.

To see the equivalence of the second formula, we compute that
$$t\frac{d}{dt}\ln\left(\frac{(-t^{-\hf};q)_\infty}{(t^{-\hf};q)_\infty}\right)
=
\frac{1}{2}\sum_{r=0}^\infty\left(\frac{q^rt^{-\hf}}{1+q^rt^{-\hf}}
+ \frac{q^rt^{-\hf}}{1-q^rt^{-\hf}}\right).$$ Combining common
denominators, one obtains the term $\frac{1}{t^\hf-t^{-\hf}}$ and
the first group of summands of the first formula.  The remaining
terms in the $\ln$ expression of the second formula yield the
second group of summands of the first formula.
\end{proof}

\subsection{The $q$-dimension of a $\dinf$-module of level $l$}
\label{sec:qdim-dinf-l}

Given an $\dinf$-module $M$, we refer to $\dim_q M := \trace_M
q^{L_0}$ as the $q$-dimension of $M$.

\begin{prop}
The $q$-dimension $\dim_q {L(d_\infty;\Lambda(\la))}$ for $\la_l
\neq 0, $ or $ \dim_q {[L(d_\infty;\Lambda(\la)) \oplus
L(d_\infty;\Lambda({\la}\otimes\det))]}$ for $\la_l = 0,$ is given
by the following (equivalent) formulas:
\begin{eqnarray*}
&& \frac1{(q;q)_\infty^{l}} \cdot
\sum_{\sigma\in W(D_l)}(-1)^{\ell(\sigma)} q^{\frac{\Vert
\la+\rho-\sigma(\rho) \Vert^2}{2}} \\
 &=& \frac1{(q;q)_\infty^{l}}  \cdot
q^{\frac{\Vert\la\Vert^2}{2}} \prod_{1\leq i<j\leq l}
\left(1-q^{\la_i-\la_j+j-i}\right)
\left(1-q^{\la_i+\la_j+2l-i-j}\right).
\end{eqnarray*}
\end{prop}

\begin{proof}
We substitute Lemma~\ref{lem:sym} with the simple identity

$$
\trace_\fock(z^{e_{11}}q^{L_0})
 =\dim_q \fock^{(0)} \sum_{k\in\Z}z^kq^{\frac{k^2}{2}},
$$
and substitute Lemma~\ref{lem:AA=D} with the identity $\dim_q
\fock^{(0)} = (q;q)_\infty^{-1}$. In this way, the same strategy
of establishing the $1$-point function in Theorem~\ref{th:1ptD}
applies and it readily leads to the first $q$-dimension formula in
the proposition. The second formula follows from
Lemma~\ref{lem:weyldenom} below and the explicit root system of
$D_l$.
\end{proof}

\begin{lemma} \label{lem:weyldenom}
Let $\mathfrak g$ be a semisimple Lie algebra with Weyl group $W$.
Set $\Vert x \Vert^2 =(x,x)$, and let $\la$ be a weight. Then,
$$\sum_{\sigma\in W}(-1)^{\ell(\sigma)} q^{\frac{\Vert
\la+\rho-\sigma(\rho) \Vert^2}{2}}
= q^{\frac{\Vert\la\Vert^2}{2}}\prod_{\alpha \in \Delta^+} (
1-q^{(\la +\rho, \alpha)}).$$
\end{lemma}
\begin{proof}
The Weyl denominator formula reads
$$\sum_{\sigma\in W}(-1)^{\ell(\sigma)} q^{\rho -\sigma(\rho)}
= \prod_{\alpha \in \Delta^+} ( 1-q^{\alpha}).$$

The lemma now follows by applying the bilinear pairing with
$\la+\rho$ to both sides of this formula of Weyl and noting that
$$\hf  {\Vert \la+\rho-\sigma(\rho) \Vert^2}
= \hf  {\Vert\la\Vert^2}  +(\la+\rho, \rho -\sigma(\rho)).
$$
\end{proof}

\section{Correlation functions on $d_\infty$-modules of level $l+\hf$}
\label{sec:dB}

\subsection{The Fock space $\fock^{l+\hf}$}
\label{sec:Fockl+hf}

Recall $\underline{\Z} = \hf + \Z \mbox{ or } \Z$. Consider the
neutral fermion
$$
\vphi(z) = \sum_{n\in \underline{\Z}}\vphi_nz^{-n-\hf+\eps}
$$
which satisfies the commutation relation
$$
[\vphi_m,\vphi_n]_+ = \delta_{m,-n}.
$$
We denote by $\fock^{l+\hf}$ the Fock space of one neutral fermion
$\vphi(z)$ and $l$ pairs of complex fermions $\psi^{\pm,p}(z),
1\le p \le l$, generated by a vacuum vector $\vac$ which satisfies
\begin{eqnarray*}
  \vphi_n \vac =\psi^{+,p}_n \vac = \psi^{-,p}_n \vac = 0
  & \text{ for } n \in \hf +\Z_+, &\mbox{ if }
   \underline{\Z} = \hf +\Z;                       \\
 \vphi_n \vac = \psi^{+,p}_n \vac = \psi^{-,p}_{n+1} \vac = 0
 &  \text{ for } n \in \Z_{+},
   & \mbox{ if } \underline{\Z} = \Z.
\end{eqnarray*}

Let
$$e^\pm_{p}  =\sum_{r\in \underline{\Z}}\no \psi_{r}^{\pm p}\vphi_{-r}\no ,
 \quad 1\le p\le l.
$$
It is known (cf. \cite{FF, W1}) that the above operators $e^+_p,
e^-_p$ together with $e^+_{pq}, e_{pq}, e^-_{pq}$ $(p, q = 1,
\cdots, l)$ defined in (\ref{e-}, \ref{e+}) generate Lie algebra
$\mathfrak {so}(2l+1)$.

\begin{lemma} \label{lem:fockiso}
Given a pair of complex fermions $\psi^{\pm} (z)$, we let
$$
\varphi_n := (\psi^+_n +\psi^-_n)/\sqrt{2},
 \qquad \varphi_n' :=  i(\psi^+_n-\psi_n^-)/\sqrt{2}.
$$
Then, $\varphi_n$ and $\varphi_n'$ satisfy the anti-commutation
relations:
\begin{eqnarray*}
 {[}\varphi_n,\varphi_m]_+ = \delta_{n,-m}, &&
 {[}\varphi_n',\varphi_m']_+ = \delta_{n,-m},\\
 {[}\varphi_n,\varphi_m'{]}_+ = 0, && \text{for } m,n\in \underline{\Z}.
\end{eqnarray*}
Hence, there is an isomorphism of Fock spaces
\[\fock^\half\otimes\fock^\half \cong \left\{
\begin{array}{ll}
\fock & \text{ if }\; \underline{\Z} =\hf +\Z, \\
\fock \oplus \fock & \text{ if }\; \underline{\Z} = \Z.
\end{array}
\right.
\]
 \end{lemma}
\begin{proof}
The commutation relations are verified by a direct computation.
The multiplicity $2$ in the Fock space isomorphism for
$\underline{\Z} = \Z$ is due to the fact that the central elements
$\varphi_0, \varphi_0'$ thus defined satisfy $\varphi_0 \vac =i
\varphi_0' \vac$.
\end{proof}
\subsection{The  $(O(2l+1), \dinf)$-Howe duality}

Now let $\underline{\Z} = \hf + \Z$ in the remainder of
Section~\ref{sec:dB}. The Lie algebra $\dinf$ acts on $\fock^l$ by
(cf. \cite{DJKM2})
\begin{eqnarray} \label{eq:genD}
&& \sum_{i,j\in\Z}(E_{i,j} -E_{1-j,1-i})z^{i-1}w^{-j} \\
 &=& \sum_{p =1}^l
\no\psi^{+,p}(z) \psi^{-,p}(w)\no -\no\psi^{+,p}(w)\psi^{-,p}(z)
\no + \no\vphi (z) \vphi (w)\no. \nonumber
\end{eqnarray}

The action of Lie algebra $\mathfrak {so}(2l+1)$ on
$\fock^{l+\hf}$ can be integrated to an action of the Lie group
$SO(2l+1)$. We identify the Borel subalgebra ${\mathfrak b}(
{\mathfrak {so}}(2l+1) )$ with the one generated by $e^+_{pq} \;(p
\neq q), e_{pq}\; (p \leq q), e^+_p,$ where $p, q = 1, \ldots, l.$
The element $\omega := \text{diag} (1, \ldots, 1, -1) \in O(2l+1)$
acts on $\fock^{l+\hf}$ by sending $\vphi_n$ to $-\vphi_n$ for
each $n$.

\begin{lemma}  \cite[Lemmas 4.2]{W1}
The action of the Lie group $O(2l+1)$ commutes with the action of
$\dinf$ on $\fock^{l+\hf}$.
\end{lemma}

Define a map $\Lambda$ from $\Sigma (B)$ to ${\dinf}_0^*$ by
sending
 $\lambda = (m_1, m_2, \ldots, m_l)$
 to
$$  \Lambda(\lambda) =
     (2l+1 - i -j) \hL^d_0 + (j-i)\hL^d_1 + \sum_{k=1}^i \hL^d_{m_k}  $$
and sending
 $\lambda = (m_1, m_2, \ldots, m_l) \bigotimes \det$
 to
$$   \Lambda(\lambda) =
    (j - i) \hL^d_0 + (2l +1-i -j)\hL^d_1 + \sum_{k =1}^i \hL^d_{m_k} $$
where
$ m_1 \geq \ldots \geq m_i > m_{i +1} = \ldots = m_j = 1 > m_{j+1}
= \ldots = m_l = 0.$

\begin{prop} \cite[Theorem 4.1]{W1}  \label{howeBd}
We have the $(O(2l+1),d_\infty)$-module decomposition:
$$
\fock^{l+\half} \cong \bigoplus_{\la\in\Sigma(B)} V_\la (O(2l+1))
\otimes L(d_\infty, \Lambda(\la)).
$$
\end{prop}

By (\ref{eq:genD}), we can write $\D(t)$ acting on $\fock^{l+\hf}$
as
$$
\D(t) = \sum_{k\in\HZ}t^k \left(\sum_{i=1}^l
(\psi^{+,i}_{-k}\psi^{-,i}_{k} + \psi^{-,i}_{-k}\psi^{+,i}_k) +
\vphi_{-k}\vphi_{k}\right).
$$

\subsection{The $n$-point $\dinf$-function of level $\hf$}
\label{sec:d+hf}

\begin{definition}
The $n$-point $\dinf$-correlation function of level $l+\hf$
associated to $\la \in \PA^l$, denoted by $\mathfrak
D^{l+\hf}_\la(q,{\bf t})$ or $\mathfrak
D^{l+\hf}_\la(q,t_1,\ldots, t_n)$, is
\begin{eqnarray*}
\trace_{L(d_\infty;\Lambda(\la)) \oplus L(d_\infty;\Lambda(\la
\otimes \det))} q^{L_0}\D(t_1)\cdots \D(t_n).
\end{eqnarray*}
(As a justification of this definition, Remark~\ref{rem:Oinf} also
applies here.)
\end{definition}

When $l=0$, $\mathfrak D^{\hf}_{(0)}(q,t_1,\ldots, t_n) =
\trace_{\fock^{\hf}} q^{L_0}\D(t_1)\cdots \D(t_n)$, by
Proposition~\ref{howeBd}. The aim of this subsection is to
determine this (unique) $n$-point $\dinf$-function of level $\hf$,
which will be used in the general level $l+\hf$ case in the
following subsection.

\begin{lemma} \label{lem:Dop}
Under the  isomorphism $\fock \cong \fock^\half\otimes\fock^\half$
in Lemma~\ref{lem:fockiso}, we have $\D(t) = \D_1(t) + \D_2(t)$,
where $\D_1(t) = \sum_{k\in\Z+\half}t^k \varphi_{-k}\varphi_k$ and
$\D_2(t) = \sum_{k\in\Z+\half} t^k\varphi_{-k}'\varphi_k'$.
\end{lemma}

\begin{proof}
A simple calculation reveals that
$$
\psi_{-k}^+\psi_k^-+\psi_{-k}^-\psi_k^+ = \varphi_{-k}\varphi_k +
\varphi_{-k}'\varphi_k'.
$$
Now the lemma follows from the definition of $\D(t)$.
\end{proof}

Given a subset $I=(i_1,\ldots, i_s) \subseteq \{1,\ldots,n\}$, we
denote by $I^c$ the complementary set to $I$, and ${\bf t}_I
=(t_{i_1},\ldots,t_{i_s})$. By convention, we let
\begin{eqnarray} \label{eq:empty}
\mathfrak D^{\hf}_{(0)}(q,{\bf t}_\emptyset) =\trace_{\fock^{\hf}}
q^{L_0} =(-q^\hf;q)_\infty
\end{eqnarray}
and recall from (\ref{eq:traceF}) that
\begin{eqnarray*}
\bff(z,q; t_1,\dots,t_n) = \trace_\fock
z^{e_{11}}q^{L_0}\D(t_1)\cdots \D(t_n).
\end{eqnarray*}
\begin{prop} \label{recursive}
We have
\begin{eqnarray} \label{eq:subset}
\bff (1,q;t_1,\ldots,t_n)
 &= \sum_{I\subseteq \{1,\ldots,n\}}
 \mathfrak D^{\hf}_{(0)}(q,{\bf t}_I)
 \mathfrak D^{\hf}_{(0)}(q,{\bf t}_{I^c}).
\end{eqnarray}

Equivalently, we have
$$
\begin{aligned}
\mathfrak D^{\hf}_{(0)} (q,{\bf t})
  &=& \hf (-q^\hf;q)_\infty^{-1}
  \left( \sum_{k\in\Z}  q^{\frac{k^2}{2}}
 \sum_{\vec{\eps}  \in \{\pm 1\}^n}  [\vec{\eps}] \cdot
 \left(\Pi{\bf t}^{\vec{\eps}} \right)^k
 F_{bo}(q;{\bf t}^{\vec{\eps}}) \right. \\
 && -\left. \sum_{\emptyset \subsetneq I \subsetneq\{1,\ldots,n\}}
 \mathfrak D^{\hf}_{(0)}(q,{\bf t}_I)
 \mathfrak D^{\hf}_{(0)}(q,{\bf t}_{I^c}) \right).
\end{aligned}
$$
\end{prop}

\begin{proof}
By Lemmas~\ref{lem:fockiso} and \ref{lem:Dop}, we have
\begin{eqnarray*}
&& \trace_\fock q^{L_0} \D(t_1) \cdots\D(t_n) \\
 &=&
\trace_{\fock^\half\otimes\fock^\half} q^{L_0}(\D_1(t_1) +
\D_2(t_1))\cdots (\D_1(t_n) + \D_2(t_n)) \\
 &=&
\sum_{\vec{i}\in\{1,2\}^n} \trace_{\fock^\half\otimes\fock^\half}
q^{L_0}\D_{i_1}(t_1)\D_{i_2}(t_2)\cdots\D_{i_n}(t_n).
\end{eqnarray*}
This is equivalent to the first formula in the theorem.

On the right-hand side of (\ref{eq:subset}), there are exactly two
terms equal to $\mathfrak D^{\hf}_{(0)} (q,t_1,\ldots,t_n)$, which
come from $I =\emptyset$ and $\{1,\ldots,n\}$.  Now the second
formula follows from (\ref{eq:empty}) and Lemma \ref{lemma:nptref}
which gives a formula for $\bff (1,q;t_1,\ldots,t_n)$.
\end{proof}

Proposition~\ref{recursive} allows for the determination, which is
recursive on $n$, of all $n$-point correlation functions
$\mathfrak D^{\hf}_{(0)}(q,t_1,\ldots,t_n)$. Note that the
$1$-point function $\mathfrak D^{\hf}_{(0)}(q,t)$ has been
computed in \cite{W2} (denoted by $S(t)$ therein) using partition
identities.
%
\begin{prop} \cite[Theorem 9]{W2}  \label{prop:w2}
The $1$-point function $\mathfrak D^{\hf}_{(0)}(q,t)$ is given by
\begin{eqnarray*}
 (-q^\hf;q)_\infty
   \left (
   \frac{1}{t^{\hf} -t^{-\hf}} +
 \sum_{r=0}^\infty
   \left[\frac{(-1)^r (q^{r+1} t)^{\hf}}{1 -q^{r+1} t}
  -\frac{(-1)^r (q^{r+1} t^{-1})^{\hf}}{1 -q^{r+1}  t^{-1}}
  \right]
   \right ).
\end{eqnarray*}
\end{prop}

An alternative solution to the $1$-point function follows from
Proposition~\ref{recursive} for $n=1$:
\begin{eqnarray}
\mathfrak D^{\hf}_{(0)}(q,t)
%
&=& (-q^\hf;q)_\infty^{-1} \sum_{k\in\Z} q^{\frac{k^2}{2}} t^k
F_{bo}(q,t) \nonumber \\
&=& \frac{(-q^\hf t;q)_\infty(-q^\hf t^{-1};q)_\infty} {
 (-q^\hf;q)_\infty \Theta(t)} \nonumber \\
&=& \frac1{(t^\hf -t^{-\hf})} \cdot \frac{(-q^\hf
t;q)_\infty(-q^\hf t^{-1};q)_\infty (q;q)_\infty^2}
{(qt;q)_\infty(qt^{-1};q)_\infty
 (-q^\hf;q)_\infty}    \label{eq:alt}
\end{eqnarray}
where we have used $F_{bo}(q,t^{-1}) =-F_{bo}(q,t)$ and the Jacobi
triple product identity. Comparing this formula with
Proposition~\ref{prop:w2} gives us the following.

\begin{cor}
The following $q$-identity holds:
$$
\begin{aligned}
& \frac{(-q^\hf t;q)_\infty(-q^\hf t^{-1};q)_\infty
(q;q)_\infty^2} {(qt;q)_\infty(qt^{-1};q)_\infty
 (-q^\hf;q)_\infty^2} \\
 &  \qquad
  = 1 + (t^\hf -t^{-\hf}) \cdot
 \sum_{r=0}^\infty
   \left[\frac{(-1)^r (q^{r+1} t)^{\hf}}{1 -q^{r+1} t}
  -\frac{(-1)^r (q^{r+1} t^{-1})^{\hf}}{1 -q^{r+1}  t^{-1}}
  \right].
 \end{aligned}
 $$
\end{cor}


\subsection{The $n$-point $\dinf$-functions of level $l+\hf$}

For $\la\in \PA^l$, the character of the irreducible
$O(2l+1)$-module associated to $\la$ and $\la\otimes \det$ is the
same, and is given as follows (cf. \cite[p. 408]{FH})
\begin{eqnarray} \label{eq:charB}
\text{ch}^b_\la(z_1,\dots,z_l)
 = \frac{\left|z_j^{\la_i +l -i +\half} -z_j^{-(\la_i +l -i
+\half)} \right|}{\left|z_j^{l-i +\half} -z_j^{-(l -i
+\half)}\right|}.
\end{eqnarray}

The following lemma is straightforward.
\begin{lemma} \label{weylcharB}
Among all the monomials $z_1^{\mu_1}\cdots z_l^{\mu_l}$ in the
expansion of the determinant ${\left|z_j^{\la_i +l -i +\half}
-z_j^{-(\la_i +l -i +\half)} \right|}$, there is exactly one
dominant monomial with $\mu_1 \ge \ldots \ge \mu_l \ge 0$, that
is, ${\bf z}^{\la+\rho}=\prod_{i=1}^l z_i^{\la_i+l-i+\hf}$. Its
coefficient is equal to $1$.
\end{lemma}

Recall the definition (\ref{eq:traceF}) of $\bff(z,q;
t_1,\dots,t_n)$.

\begin{lemma}  \label{lem:howeBd}
We have the following $q$-series identity:
\begin{eqnarray*}
&&  \trace_{\fock^{\half}}q^{L_0} \D(t_1)\cdots\D(t_n) \cdot
\prod_{i=1}^l  \bff(z_i,q; t_1,\dots,t_n)
\\
&=& \sum_{\la \in\PA^l} \text{ch}^b_\la(z_1,\dots,z_l) \mathfrak
D_\la^{l+\half}(q;t_1,\ldots,t_n).
\end{eqnarray*}
\end{lemma}

\begin{proof}
This follows from the application of
 $
\trace_{\fock^{l+\half}} z_1^{e_{11}} \cdots z_l^{e_{ll}} q^{L_0}
\D(t)
 $
to both sides of the Howe duality in Proposition~\ref{howeBd}.
Note that $\fock^{l+\hf} \cong \fock^{l} \otimes \fock^{\half}$.
On the left-hand side, $z_i^{e_{ii}}$ only acts on the $i^{th}$
tensor factor of $\fock^{l}$ and not on $\fock^{\half}$. For the
right-hand side, $z_1^{e_{11}}\cdots z_l^{e_{ll}}$ only acts on
the first tensor factor and $q^{L_0}\D(t)$ only acts on the second
tensor factor.
\end{proof}

Recall that $\mathfrak D^{\hf}_{(0)} (q;{\bf t})$ has been
computed recursively in the previous subsection.
\begin{theorem} \label{th:lplushalf}
The $n$-point $\dinf$-correlation function of level $l+\hf$,
$\mathfrak D^{l+\hf}_\la(q,t_1,\ldots,t_n)$, is equal to
\begin{eqnarray*}
 && \mathfrak D^{\hf}_{(0)} (q;{\bf t}) \times \\
 &&\times
 \sum_{\sigma\in W(B_l)}\left(-1\right)^{\ell(\sigma)}
q^{\frac{\Vert \la+\rho-\sigma(\rho) \Vert^2}{2}}
\prod_{a=1}^l
 \Big(
\sum_{\vec{\eps}_a\in\{\pm 1\}^n}
 [\vec{\eps}_a] (\Pi {\bf t}^{\vec{\eps}_a})^{k_a}
 F_{bo}(q;{\bf t}^{{\vec{\eps}_a}})
 \Big)
\end{eqnarray*}
where $k_a = (\la+\rho-\sigma(\rho), \ep_a)$.
\end{theorem}

\begin{proof}
 The Weyl denominator of type $B_l$ reads that
\begin{eqnarray} \label{eq:denomB}
 {\left |z_j^{l-i+\hf}+z_j^{-(l-i+\hf)}\right|} = \sum_{\sigma\in
W(B_l)} (-1)^{\ell(\sigma)} {\bf z}^{\sigma(\rho)}.
\end{eqnarray}

It follows by (\ref{eq:charB}), (\ref{eq:denomB}),
Lemmas~\ref{lemma:nptref} and \ref{lem:howeBd} that
\begin{eqnarray*}
&& \sum_{\sigma\in W(B_l)} (-1)^{\ell(\sigma)} {\bf
z}^{\sigma(\rho)} \cdot \mathfrak D_{(0)}^\hf(q;{\bf t}) \times \\
 && \quad \qquad \times
 \prod_{a=1}^l
 \left( \sum_{k_a\in\Z} z_a^{k_a} q^{\frac{k_a^2}{2}}
 \sum_{\vec{\eps}_a  \in \{\pm 1\}^n}  [\vec{\eps}_a] \cdot
 \left(\Pi{\bf t}^{\vec{\eps}_a} \right)^{k_a}
 F_{bo}(q;{\bf t}^{{\vec{\eps}_a}})
 \right)
 \nonumber \\
 &=& \sum_{\la \in \PA^l}  {\left|z_j^{\la_i+l-i+\hf} +
z_j^{-(\la_i+l-i+\hf)}\right|} \cdot \mathfrak D^{l+\hf}_\la(q,{\bf
t}).
\end{eqnarray*}

Now the theorem follows by Lemma~\ref{weylcharB} and by comparing
the coefficients of ${\bf z^{\la+\rho}}$ on both sides of the
above identity.
\end{proof}
\subsection{The $q$-dimension of a $\dinf$-module of level $l+\hf$}

In the same manner as in Section~\ref{sec:qdim-dinf-l}, we can
derive the following $q$-dimension formula from the $(O(2l+1),
\dinf)$-Howe duality in Proposition~\ref{howeBd}. The second
formula below is obtained from the first one by using
Lemma~\ref{lem:weyldenom} and the explicit root system of type
$B_l$. Recall that $L(d_\infty;\Lambda(\la)) \oplus
L(d_\infty;\Lambda({\la}\otimes\det))$ can be regarded as an
irreducible module of the orthogonal group corresponding to
$\dinf$ (cf. Remark~\ref{rem:Oinf}).

\begin{prop}
We have
\begin{align*}
 \dim_q & [L(d_\infty;\Lambda(\la)) \oplus
L(d_\infty;\Lambda({\la}\otimes\det))]  \\
 & =
\frac{(- q^{-\hf};q)_\infty}{(q;q)_\infty^l} \cdot
\sum_{\sigma\in W(B_l)}(-1)^{\ell(\sigma)} q^{\frac{\Vert
\la+\rho-\sigma(\rho) \Vert^2}{2}} \\
 &=
\frac{(- q^{-\hf};q)_\infty}{(q;q)_\infty^l} \cdot
q^{\frac{\Vert\la\Vert^2}{2}}
  \prod_{1\leq i\leq l} \left(1-q^{\la_i+l-i+1/2}\right) \times\\
 & \quad \times \prod_{1\leq i<j\leq l}
\left(1-q^{\la_i-\la_j+j-i}\right)
\left(1-q^{\la_i+\la_j+2l-i-j+1}\right).
\end{align*}
\end{prop}

\section{Correlation functions on $\cinf$-modules of level $l$}
\label{sec:cC}

\subsection{The  $(Sp(2l),\cinf)$-Howe duality}

We again take $\underline{\Z} = \hf+\Z$ for the Fock space
$\fock^l$ of fermions $\psi^{\pm,p}(z)$, $1\le p \le l$, in this
section. The representation of $\cinf$ on $\fock^l$ is given by
(\cite{DJKM2})
\begin{eqnarray}\label{gs:cinf}
&& \sum_{i,j\in\Z}(E_{i,j}-(-1)^{i+j} E_{1-j,1-i})z^{i-1}w^{-j}
 \nonumber  \\
 &=& \sum_{p =1}^l
\no\psi^{+,p}(z) \psi^{-,p}(w)\no -\no\psi^{+,p}(w)\psi^{-,p}(z)
\no.
\end{eqnarray}

Let
$$\tilde{e}^-_{pq}
=\sum_{r\in\hf+ \Z}(-1)^{r-\hf}\no\psi_{r}^{-p}\psi^{-q}_{-r}\no,\quad
\tilde{e}_{pq}^+ = \sum_{r\in\hf+
\Z}(-1)^{r-\hf}\no\psi_{r}^{+p}\psi^{+q}_{-r}\no,$$ and let
$$
\tilde{e}_{pq}  =\sum_{r\in\hf+
\Z}\no\psi_{r}^{+p}\psi^{-q}_{-r}\no.
$$
The operators $\tilde{e}^+_{pq}, \tilde{e}_{pq}, \tilde{e}^-_{pq}
\;(p, q = 1, \ldots, l)$ generate Lie algebra $\mathfrak {sp}(2l)$
and can be integrated to the action of the Lie group $Sp(2l)$ on
$\fock^l$ (cf. \cite{FF, W1}). In particular, the operators
$\tilde{e}_{pq}\; ( p, q = 1, \ldots, l)$ form a Lie subalgebra
${\mathfrak {gl}}(l)$ in the horizontal of $\mathfrak {sp}(2l)$.
Identify the Borel subalgebra ${\mathfrak b}( {\mathfrak {sp}}(2l)
)$ with the one generated by $\tilde{e}_{pq}\;(p \leq q),
\tilde{e}_{pq}^+, p, q = 1, \ldots, l.$ It is known by
\cite[Lemmas 3.6, Remark 3.7]{W1} that the action of the Lie group
$Sp(2l)$ commutes with the action of $\cinf$ on $\fock^l$.

Define the map $\Lambda: \PA^l \longrightarrow {\cinf}_0^* $ by
sending $\lambda = (m_1, \ldots, m_l)$ to
$ \Lambda (\lambda) =
   (l-j) \hL_0^c + \sum_{k =1}^j \hL_{m_k}^c, $
where $j$ denotes the last non-zero index among $m_i$'s.

\begin{prop} \cite[Theorem 3.4]{W1} \label{prop:HoweCc}
We have the following decomposition of
$(Sp(2l),c_\infty)$-modules:
\begin{eqnarray}\label{c:duality}
\fock^{l} = \bigoplus_{\la\in\PA^l} V_\la(Sp(2l))\otimes
L(c_\infty;\Lambda(\la)).
\end{eqnarray}
\end{prop}

\subsection{The $n$-point $\cinf$-functions of level $l$}

Introduce the following operators in $\cinf$:
\begin{eqnarray*}
 \no\CC(t)\no &=& \sum_{k\in \mathbb N}(t^{k-\hf} -
t^{\hf-k})(E_{k,k}-E_{1-k,1-k}), \\
 \CC(t) &=& \no \CC(t)\no +\frac{2}{t^\hf-t^{-\hf}}C.
\end{eqnarray*}
When acting on $\fock^l$, these operators can be written in terms
of the operators $\psi^\pm_n$ by (\ref{gs:cinf}) as
\begin{eqnarray*}
\no \CC(t)\no &=&\sum_{p =1}^l \sum_{k\in\HZ} t^k (\no
\psi^{+,p}_{-k} \psi^{-,p}_{k}\no + \no\psi^{-,p}_{-k}
\psi^{+,p}_k \no), \\
 \CC(t) &=& \sum_{p =1}^l \sum_{k\in\HZ} t^k
(\psi^{+,p}_{-k} \psi^{-,p}_{k} +\psi^{-,p}_{-k} \psi^{+,p}_k).
\end{eqnarray*}

\begin{definition}
The $n$-point $\cinf$-correlation function of level $l$ associated
to $\la  \in \PA^l$ is
$$
\mathfrak C^l_\la(q,t_1,\ldots, t_n) =
\trace_{L(c_\infty;\Lambda(\la))} q^{L_0}\CC(t_1)\cdots \CC(t_n).
$$
\end{definition}

\begin{theorem}
The $n$-point $\cinf$-correlation function of level $l$ is given
by
\begin{eqnarray*}
&& \mathfrak C_\la^l(q;t_1,\dots,t_n) =\\
&& \sum_{\sigma\in W(C_l)}\left(-1\right)^{\ell(\sigma)}
q^{\frac{\Vert \la+\rho-\sigma(\rho) \Vert^2}{2}}
 \prod_{a=1}^l
 \Big(
\sum_{\vec{\eps}_a\in\{\pm 1\}^n}
 [\vec{\eps}_a] (\Pi {\bf t}^{\vec{\eps}_a})^{k_a}
 F_{bo}(q;{\bf t}^{{\vec{\eps}_a}})
 \Big)
\end{eqnarray*}
where $k_a = (\la+\rho-\sigma(\rho), \ep_a)$.
\end{theorem}

\begin{proof}
Note that $\CC(t) = \A(t) - \A(t^{-1})$. The proof follows the same
strategy which works for Theorems~\ref{th:1ptD} and \ref{th:nlD} for
$\dinf$-correlation functions of level $l$. We now use instead the
combinatorial consequence of the $(Sp(2l),c_\infty)$-Howe duality
(\ref{c:duality}) and the character of irreducible $Sp(2l)$-modules
(cf. \cite[24.18]{FH})
$$
\text{ch}_\la^{sp}(z_1,\dots,z_l) = \frac{\left| z_j^{\la_i+l-i+1}
-z_j^{-(\la_i+l-i+1)}\right |}{\left |z_j^{l-i+1}-z_j^{-(l-i+1)}
\right |}.
$$
Note that the Weyl group $W(C_l)$ replaces $W(D_l)$ in the proof
and result.
\end{proof}

In the case $n=1$, the notation can be much simplified. The
$1$-point $\cinf$-function of level $l$ is given by
$$
\mathfrak C^l_\la(q,t) = F_{bo}(q,t)^l \cdot
\sum_{\sigma\in W(C_l)}(-1)^{\ell(\sigma)} q^{\frac{\Vert
\la+\rho-\sigma(\rho) \Vert^2}{2}}
\prod_{a=1}^l(t^{k_a}+t^{-k_a})$$
where $k_a = (\la+\rho-\sigma(\rho), \ep_a)$.

Let us specialize further to $l=1$.  The irreducible character of
$Sp(2) =SL(2)$ is simply
 $$
\text{ch}_m^{sp}(z) = (z^{m+1}-z^{-(m+1)})/ (z-z^{-1}).
 $$
Then the 1-point $\cinf$-correlation function of level $1$ is
given by
$$
\mathfrak C^1_{(m)}(q,t) = \frac{q^{m^2/2} \left(t^m+t^{-m}\right)
-q^{(m+2)^2/2}\left(t^{m+2}
+t^{-(m+2)}\right)}{(q;q)_\infty\Theta(t)}.
$$
In contrast to the $\dinf$ case at level $1$ where the charge
decomposition of $\fock$ and the theory of partitions can be used
effectively, the description of irreducible $\cinf$-submodules in
$\fock$ is not explicit and the Howe duality in
Proposition~\ref{prop:HoweCc} is essentially used.

\subsection{The $q$-dimension of a $\cinf$-module}

In the same manner as in Section~\ref{sec:qdim-dinf-l}, we can
derive the following $q$-dimension formula from the $(Sp(2l),
\cinf)$-Howe duality in Proposition~\ref{prop:HoweCc}. The second
formula below is obtained from the first one by using
Lemma~\ref{lem:weyldenom} and the explicit root system of type
$C_l$.

\begin{prop}
For $\la \in \PA^l$, we have
\begin{align*}
\dim_q L (\cinf;\Lambda(\la))
 & =
\frac{1}{(q;q)_\infty^l} \cdot
\sum_{\sigma\in W(C_l)}(-1)^{\ell(\sigma)} q^{\frac{\Vert
\la+\rho-\sigma(\rho) \Vert^2}{2}} \\
 &=
\frac1{(q;q)_\infty^l} \cdot
q^{\frac{\Vert\la\Vert^2}{2}}
  \prod_{1\leq i\leq l} \left(1-q^{2(\la_i+l-i+1)}\right) \times\\
 &\quad \times \prod_{1\leq i<j\leq l}
\left(1-q^{\la_i-\la_j+j-i}\right)
\left(1-q^{\la_i+\la_j+2l-i-j+2}\right).
\end{align*}
\end{prop}

\section{Correlation functions on $\binf$-modules of level $l$}
\label{sec:binfl}
\subsection{The $(Pin(2l),\binf)$-Howe duality}

Throughout Section~\ref{sec:binfl} we take $\underline{\Z} = \Z$.
The action of $\binf$ on the Fock space $\fock^l$ is given by
(\cite{DJKM2})
\begin{equation}\label{gs:binf}
\sum_{i,j\in\Z}(E_{i,j}-E_{-j,-i})z^{i}w^{-j} = \sum_{p=1}^l
 \left( \no\psi^{+,p}(z)\psi^{-,p}(w)\no
-\no\psi^{+,p}(w)\psi^{-,p}(z)\no
 \right).
\end{equation}
The Lie algebra $\mathfrak{so}(2l)$ defined in
Section~\ref{sec:HoweDd} can be integrated to $Spin(2l)$ and then
naturally extended to $Pin(2l)$. Remark 3.5 and Lemma 3.5 of
\cite{W1} are summed up by the following lemma.
\begin{lemma}
The action of the Lie group $Pin(2l)$ commutes with the action of
$\binf$ on $\fock^l$.
\end{lemma}
For $\la = {\bf 1}_l/2 + (m_1,\dots,{m}_l)$ in $\Sigma(Pin)$,
define the following map $\Lambda: \Sigma(Pin)\longrightarrow
{\binf}_0^*$:
$$\Lambda(\la) = (2l-2j)\hL_0^b + \sum_{k=1}^j \hL_{m_k}^b$$
where $j$ is such that $m_1\geq\cdots\geq m_j > m_{j+1}
=\cdots=m_l = 0$.
\begin{prop} \cite[Theorem 3.3]{W1} \label{prop:HoweDb}
We have the following decomposition of $(Pin(2l),\binf)$-modules:
\begin{eqnarray*}
\fock^{l} = \bigoplus_{\la\in\Sigma(Pin)} V_\la(Pin(2l))\otimes
L(\binf;\Lambda(\la)).
\end{eqnarray*}
\end{prop}
\subsection{The operator $\B(t)$}

Introduce the following operators in $\binf$:

\begin{eqnarray*}
 \no\B(t)\no &=& \sum_{k\in \Z_+}(t^k-t^{-k})(E_{k,k}-E_{-k,-k}), \\
 \B(t) &=& \no\B(t)\no +\frac{t+1}{t-1}C.
\end{eqnarray*}
When acting on $\fock^l$, $\no\B(t)\no$ and $\B(t)$ can be
expressed using (\ref{gs:binf}) as follows:
\begin{eqnarray*}
\no\B(t)\no &=&\sum_{p=1}^l \sum_{k\in\Z} t^k (\no\psi^{+,p}_{-k}
\psi^{-,p}_{k}\no + \no\psi^{-,p}_{-k}
\psi^{+,p}_k\no), \\
 \B(t) &=& \sum_{p=1}^l \sum_{k\in\Z} t^k
(\psi^{+,p}_{-k} \psi^{-,p}_{k} +\psi^{-,p}_{-k} \psi^{+,p}_k).
\end{eqnarray*}
We easily verify that
\begin{equation} \label{eq:AB}
\B(t) = t^\hf\A(t) - t^{-\hf} \A(t^{-1}), \qquad \no \B(t)\no =
t^\hf \no \A(t)\no - t^{-\hf} \no \A(t^{-1})\no.
\end{equation}

The energy operator $L_0$ on the $\binf$-module
$L(\binf,\hL(\la))$ with highest weight vector $v_{\hL(\la)}$ is
defined by (\ref{weight}) and
$$
 {[}L_0, E_{i,j}-E_{-j,-i}] = (i-j) (E_{i,j}-E_{-j,-i}).
$$
On $\fock^l$, we can realize $L_0$ as
\begin{eqnarray}  \label{eq:L0binf}
 L_0 =\sum_{p=1}^l \sum_{k\in\Z} k
 \no\psi_{-k}^{+,p}\psi_k^{-,p}\no +\frac{l}{8}.
\end{eqnarray}

\subsection{The $n$-point $\binf$-functions
of level $l$}

\begin{definition}
The $n$-point $\binf$-correlation function of level $l$ associated
to $\la =\mathbf{1}_l/2 + (m_1, \ldots, m_l) \in \Sigma (Pin)$ is
$$\mathfrak B^l_\la(q,{\bf t}\,) \equiv \mathfrak
B^l_\la(q,t_1,\ldots, t_n) = \trace_{L(\binf;\Lambda(\la))}
q^{L_0}\B(t_1)\cdots \B(t_n).$$
\end{definition}

\begin{lemma} \label{lem:charspin2l}
Let $\la \in \Sigma(Pin)$ and denote by $\text{ch}^{pin}_\la
(z_1,\dots,z_l)$ the character of $V_\la(Pin(2l))$. Then,
$$
\text{ch}^{pin}_\la (z_1,\dots,z_l)
 = \frac{2\left|z_j^{\la_i+l-i} +
z_j^{-(\la_i+l-i)} \right|}{\left| z_j^{l-i}+z_j^{-(l-i)}\right|}.
$$
\end{lemma}

\begin{proof}
Recall that $\overline{\la} = \overline{\bf 1}_l/2 +
(m_1,\dots,m_{l-1},-{m}_l)$. Since $V_\la(Pin(2l)) \cong
V_\la(Spin(2l)) \oplus V_{\overline{\la}}(Spin(2l))$, we have $
\text{ch}^{pin}_\la = \text{ch}^{so}_\la +
\text{ch}^{so}_{\overline{\la}}$, where
 $\text{ch}^{so}_\la$ for $\la \in \Sigma(Pin)$ is also
given by (\ref{char:so}). Note that the second determinant terms
in the numerators of $\text{ch}^{so}_\la$ and
$\text{ch}^{so}_{\overline{\la}}$ (cf. (\ref{char:so})) are
opposite to each other. Now the formula for $\text{ch}^{pin}_\la$
follows.
\end{proof}
Let
\begin{eqnarray} \label{eq:traceFB}
\bff_b(z,q; t_1,\dots,t_n) := \trace_\fock
z^{e_{11}}q^{L_0}\B(t_1)\cdots \B(t_n).
\end{eqnarray}

\begin{lemma} \label{lem:nptB}
We have
\begin{eqnarray*}
 \bff_b (z,q; t_1,\dots,t_n) = \sum_{k\in\hf +\Z} z^k q^{\frac{k^2}{2}}
 \sum_{\vec{\eps}  \in \{\pm 1\}^n}  [\vec{\eps}] \cdot
 \left(\Pi{\bf t}^{\vec{\eps}} \right)^k
 F_{bo}(q;{\bf t}^{\vec{\eps}})
\end{eqnarray*}
where the notations are as in Lemma~\ref{lemma:nptref}.
\end{lemma}
\begin{proof}
The proof is similar to the one for Lemma~\ref{lemma:nptref}, while
we have to take into account the difference coming from the {\em
integral} indices on the fermions $\psi^{\pm}(z)$. Denote the charge
operator $C =\sum_{k\in\Z} \no\psi^+_{-k} \psi^-_k\no$ which acts as
$0$ on $\fock^{(0)}$. By definition, $e_{11} =C +\hf$. We can check
(compare \cite{MJD} and \cite[Appendix A]{Ok}) that
$$
\mathsf S^{-k} e_{11} \mathsf S^k = e_{11}+k, \quad
 \mathsf S^{-k} \A(t) \mathsf S^k = t^k\A(t),
$$
and
$$\mathsf S^{-k} L_0 \mathsf S =L_0 +k C +\hf k(k+1) =L_0 +k e_{11}
+\hf k^2.
$$
Then, using (\ref{eq:traceFB}) and (\ref{eq:AB}), we have
$$
\begin{aligned}
& \bff_b(z,q;t_1,\dots,t_n) \\
%
&= \sum_{\vec{\eps}\in\{\pm 1\}^n} [\vec{\eps}]
 \left(\Pi {\bf t}^{ \vec{\eps}}\right)^\hf
 \sum_{k\in\Z} \trace_{\fock^{(k)}} z^{e_{11}}q^{L_0}
\A(t_1^{\eps_1})\cdots\A(t_n^{\eps_n}) \\
&= \sum_{\vec{\eps}\in\{\pm 1\}^n} [\vec{\eps}]
 \left(\Pi {\bf t}^{ \vec{\eps}}\right)^\hf
 \sum_{k\in\Z} \trace_{\fock^{(0)}} \mathsf S^{-k} z^{e_{11}}q^{L_0}
\A(t_1^{\eps_1})\cdots\A(t_n^{\eps_n})\mathsf S^k \\
&=\sum_{\vec{\eps}\in\{\pm 1\}^n} [\vec{\eps}]
 \left(\Pi {\bf t}^{ \vec{\eps}}\right)^\hf
 \sum_{k\in\Z} z^{k+\hf} q^{\hf k(k+1)}
 \trace_{\fock^{(0)}}q^{L_0}\left(\Pi
 {\bf t}^{\vec{\eps}}\right)^k\A(t_1^{\eps_1})\cdots\A(t_n^{\eps_n}) \\
&= \sum_{z\in\Z} z^{k+\hf} q^{\hf (k+\hf)^2} \sum_{\vec{\eps} \in
\{\pm 1\}^n} [\vec{\eps}] \cdot
 \left(\Pi{\bf t}^{\vec{\eps}} \right)^{k+\hf}
 F_{bo}(q;{\bf t}^{\vec{\eps}}).
\end{aligned}
$$
In the last equation we have used the ``correction term" in
(\ref{eq:L0binf}).
\end{proof}

\begin{lemma} \label{lem:qbinf}
We have the following identity:
$$\prod_{i=1}^l \bff_b (z_i,q;t_1,\ldots,t_n)
 =  \sum_{\la \in \Sigma(Pin)}
  \text{ch}^{pin}_\la (z_1,\dots,z_l)\cdot \mathfrak
B^l_\la(q,t_1,\ldots,t_n).$$
\end{lemma}

\begin{proof}
Follows from the $(Pin(2l),\binf)$-Howe duality in
Proposition~\ref{prop:HoweDb}.
\end{proof}

\begin{theorem} \label{th:nlB}
The $n$-point $\binf$-correlation function of level $l$ associated
to $\la \in \Sigma (Pin)$ is given by
\begin{eqnarray*}
&& \mathfrak B_\la^l(q;t_1,\dots,t_n) =\\
&& \sum_{\sigma\in W(D_l)}\left(-1\right)^{\ell(\sigma)}
q^{\frac{\Vert \la+\rho-\sigma(\rho) \Vert^2}{2}}
 \prod_{a=1}^l
 \Big(
\sum_{\vec{\eps}_a\in\{\pm 1\}^n}
 [\vec{\eps}_a] (\Pi {\bf t}^{\vec{\eps}_a})^{k_a}
 F_{bo}(q;{\bf t}^{{\vec{\eps}_a}})
 \Big)
\end{eqnarray*}
where $k_a = (\la+\rho-\sigma(\rho), \ep_a)$.
\end{theorem}

\begin{proof}
Note that Lemma~\ref{weylcharindep} on the dominant monomial of
the numerator of $\text{ch}^{pin}_\la$ remains valid for $\la =
{\bf 1}_l/2 + (m_1,\dots,{m}_l)$. Now the proof of the theorem is
the same as for Theorems~\ref{th:1ptD} and \ref{th:nlD}, with the
help of Lemmas~\ref{lem:charspin2l}, \ref{lem:nptB}, and
\ref{lem:qbinf}.
\end{proof}

\begin{remark}
It is remarkable that the formula for $\mathfrak B_\la^l$ in
Theorem~\ref{th:nlB} coincides with the one for $\mathfrak
D_\la^l$ given in Theorem~\ref{th:nlD}, except that the $\la$ used
in these two cases are different.
\end{remark}

\subsection{The $q$-dimension of a $\binf$-module}

In the same manner as in Section~\ref{sec:qdim-dinf-l}, we can
derive the following $q$-dimension formula from the $(Pin(2l),
\binf)$-Howe duality in Proposition~\ref{prop:HoweCc}.

\begin{prop}
For $\la  \in \Sigma (Pin)$, we have
\begin{align*}
  \dim_q L &(\binf;\Lambda(\la)) \\
 & =
\frac{1}{(q;q)_\infty^l} \cdot
\sum_{\sigma\in W(D_l)}(-1)^{\ell(\sigma)} q^{\frac{\Vert
\la+\rho-\sigma(\rho) \Vert^2}{2}} \\
 &=
\frac1{(q;q)_\infty^l} \cdot
q^{\frac{\Vert\la\Vert^2}{2}} \prod_{1\leq i<j\leq l}
\left(1-q^{\la_i-\la_j+j-i}\right)
\left(1-q^{\la_i+\la_j+2l-i-j}\right).
\end{align*}
\end{prop}

\section{Correlation functions on $\binf$-modules of level $l+\half$}
\label{sec:bB}
\subsection{The $(Spin(2l+1),\binf)$-Howe duality}

Let $\underline{\Z} =\Z$ throughout this section. The action of
Lie algebra $\mathfrak{so}(2l+1)$ on the Fock space
$\fock^{l+\hf}$ defined in Section~\ref{sec:Fockl+hf} can be
integrated to an action of Lie group $Spin(2l+1)$. The action of
$\binf$ on the Fock space $\fock^{l+\hf}$ is given by
\begin{eqnarray}\label{gs:binf+hf}
&&\sum_{i,j\in\Z}(E_{i,j}-E_{-j,-i})z^{i}w^{-j} \nonumber  \\
&=& \sum_{p=1}^l \no\psi^{+,p}(z)\psi^{-,p}(w)\no
-\no\psi^{+,p}(w)\psi^{-,p}(z)\no +\no \vphi(z)\vphi(w)\no.
\end{eqnarray}
Now the operator $\B(t)$ acting on $\fock^{l+\hf}$ can be written
as
$$\B(t) = \sum_{p=1}^l \sum_{k\in\Z}
t^k\left(\psi_{-k}^{+,p}\psi_k^{-,p}
 + \psi_{-k}^{-,p} \psi_k^{+,p} \right)
+ t^k \vphi_{-k}\vphi_k.$$

Define $\Lambda': \Sigma(Pin) \rightarrow {\binf}_0^*$ by sending
$\lambda = \hf {\bf 1}_l + (m_1, m_2, \ldots, m_l)$ to
$$
\Lambda' (\lambda)
 = (2l +1 - 2j ) \hL^b_0 + \sum_{k=1}^j \hL^b_{m_k}
$$
if $m_1 \geq \cdots \geq m_j > m_{j+1} = \cdots = m_l =0.$

\begin{prop} \cite[Theorem 4.2]{W1} \label{howeBb}
We have the $(Spin(2l+1),\binf)$-module decomposition:
$$
\fock^{l+\half} \cong 2 \bigoplus_{\la\in\Sigma(Pin)}
V_\la(Spin(2l+1)) \otimes L(\binf, \Lambda'(\la))
$$
where the factor $2$ denotes the multiplicity.
\end{prop}

The energy operator $L_0$ on the $\binf$-module
$L(\binf,\hL'(\la))$ with highest weight vector $v_{\hL'(\la)}$ is
defined by
\begin{eqnarray*}
L_0 \cdot v_{\hL '(\la)}
 &=& \left(\hf \Vert\la\Vert^2 +\frac1{16}\right) \cdot v_{\hL '(\la)},  \\
{[}L_0, E_{i,j}-E_{-j,-i}] &=& (i-j) (E_{i,j}-E_{-j,-i}).
\end{eqnarray*}
The convention of shift by $\frac1{16}$ will be convenient later
on, and it also fits with the standard realization in terms of
neutral fermions with integral indices (i.e. Ramond sector) of
$L_0$ of the Virasoro algebra.

\subsection{The $n$-point $\binf$-function of level $\hf$}

\begin{definition}
The $n$-point $\binf$-correlation function of level $l+\hf$
associated to $\la \in \Sigma(Pin)$ is
\begin{eqnarray*}
\mathfrak B^{l+\hf}_\la(q,t_1,\ldots, t_n)
=\trace_{L(b_\infty;\Lambda'(\la))} q^{L_0}\B(t_1)\cdots \B(t_n).
\end{eqnarray*}
\end{definition}
On $\fock^{l+\hf}$, we can realize $L_0$ as
\begin{eqnarray}  \label{eq:L0binf-hf}
 L_0 =\sum_{p=1}^l \sum_{k\in\Z} k
 \no\psi_{-k}^{+,p}\psi_k^{-,p}\no + \sum_{k\in\Z} \frac{k}2
 \no\vphi_{-k} \vphi_k \no +\frac{2l+1}{16}.
\end{eqnarray}

When $l=0$, we have by Proposition~\ref{howeBb} that
\begin{eqnarray}  \label{eq:Bhf}
\mathfrak B^{\hf}_{(\hf)}(q,t_1,\ldots, t_n) = \hf
\trace_{\fock^{\hf}} q^{L_0}\B(t_1)\cdots \B(t_n).
\end{eqnarray}
The aim of this subsection is to determine this (unique) $n$-point
$\binf$-function of level $\hf$ parallel to
Section~\ref{sec:d+hf}, which will be used in the general level
$l+\hf$ case in the following subsection.

The following lemma is straightforward once we recall the setup of
Lemma~\ref{lem:fockiso}.

\begin{lemma} \label{lem:Bop}
Under the isomorphism
$2 \fock \cong \fock^\half\otimes\fock^\half$ in
Lemma~\ref{lem:fockiso}, we have $\B(t) = \B_1(t) + \B_2(t)$, where
$\B_1(t) = \sum_{k\in\Z}t^k \varphi_{-k}\varphi_k$ and $\B_2(t) =
\sum_{k\in\Z} t^k\varphi_{-k}'\varphi_k'$.
\end{lemma}

By convention, we let
\begin{eqnarray} \label{eq:zero}
\mathfrak B^{\hf}_{(\hf)}(q,{\bf t}_\emptyset) = \hf
\trace_{\fock^{\hf}} q^{L_0} =q^{\frac1{16}} (-q;q)_\infty.
\end{eqnarray}

\begin{prop} \label{Brecursive}
Recalling $\bff_b$ from (\ref{eq:traceFB}), we have
\begin{eqnarray} \label{eq:Bsubset}
\bff_b (1,q;t_1,\ldots,t_n)
 = 2 \sum_{I\subseteq \{1,\ldots,n\}}
 \mathfrak B^{\hf}_{(\hf)}(q,{\bf t}_I)
 \mathfrak B^{\hf}_{(\hf)}(q,{\bf t}_{I^c}).
\end{eqnarray}

Equivalently, $\mathfrak B^{\hf}_{(\hf)} (q,{\bf t})$ is given by
$$
\begin{aligned}
&& \hf q^{-\frac1{16}}(-q;q)_\infty^{-1}
  \left(\hf \sum_{k\in\Z+\hf}  q^{\frac{k^2}{2}}
 \sum_{\vec{\eps}  \in \{\pm 1\}^n}  [\vec{\eps}] \cdot
 \left(\Pi{\bf t}^{\vec{\eps}} \right)^k
 F_{bo}(q;{\bf t}^{\vec{\eps}}) \right. \\
 && -\left. \sum_{\emptyset \subsetneq I \subsetneq\{1,\ldots,n\}}
 \mathfrak B^{\hf}_{(\hf)}(q,{\bf t}_I)
 \mathfrak B^{\hf}_{(\hf)}(q,{\bf t}_{I^c}) \right).
\end{aligned}
$$
\end{prop}

\begin{proof}
By Lemmas~\ref{lem:fockiso} and \ref{lem:Bop}, we have
\begin{align*}
 2 \trace_\fock & q^{L_0} \B(t_1) \cdots\B(t_n) \\
 &=
\trace_{\fock^\half\otimes\fock^\half} q^{L_0}(\B_1(t_1) +
\B_2(t_1))\cdots (\B_1(t_n) + \B_2(t_n)) \\
 &=
\sum_{\vec{i}\in\{1,2\}^n} \trace_{\fock^\half\otimes\fock^\half}
q^{L_0}\B_{i_1}(t_1)\B_{i_2}(t_2)\cdots\B_{i_n}(t_n).
\end{align*}
This is equivalent to the first formula in the theorem by
(\ref{eq:Bhf}), since the definitions of $L_0$ on $\fock$ and on
${\fock^\hf \otimes \fock^\hf}$ are compatible (see
(\ref{eq:L0binf}) and (\ref{eq:L0binf-hf})).

On the right-hand side of (\ref{eq:Bsubset}), there are exactly
two terms from $I =\emptyset$ and $\{1,\ldots,n\}$ which give rise
to $\mathfrak B^{\hf}_{(\hf)} (q,t_1,\ldots,t_n)$.  Now the second
formula follows from (\ref{eq:zero}) and Lemma \ref{lem:nptB}.
\end{proof}

Proposition~\ref{Brecursive} allows for the determination, which
is recursive on $n$, of all $n$-point correlation functions
$\mathfrak B^{\hf}_{(\hf)}(q,t_1,\ldots,t_n)$ of level $\hf$. The
$1$-point function $\mathfrak B^{\hf}_{(\hf)}(q,t)$ has been
computed in \cite{W2} (denoted by $R(t)$ therein up to a factor
$q^{\frac1{16}}$) using partition identities.

\begin{prop} \cite[Theorem 4]{W2} \label{w2:b}
The $1$-point function $\mathfrak B^{\hf}_{(\hf)}(q,t)$ is given
by
\begin{eqnarray*}
 q^{\frac1{16}}  (-q;q)_\infty
 \left(
 \frac{t+1}{2(t-1)} +\sum_{r =0}^\infty
 \left[ \frac{(-1)^r q^{r+1} t}{1 -q^{r+1} t}
  -\frac{(-1)^r q^{r+1} t^{-1}}{1 -q^{r+1} t^{-1}}
 \right]
 \right ).
\end{eqnarray*}
\end{prop}

An alternative solution to the 1-point function follows from
Proposition~\ref{Brecursive} for $n=1$:
$$
\begin{aligned}
\mathfrak B^{\hf}_{(\hf)}(q,t)
 &=
\frac{1}{4q^{\frac{1}{16}}(-q;q)_\infty} \sum_{k\in\Z+\hf}
q^{\frac{k^2}{2}} (t^k+t^{-k})F_{bo}(q,t)
 \\ &=
 \frac{q^{\frac18}t^\hf
(-qt;q)_\infty(-t^{-1};q)_\infty (q;q)_\infty}{2
q^{\frac{1}{16}}(-q;q)_\infty \cdot t^\hf(1-t^{-1})
\Theta(t)(q;q)_\infty}
 \\ &=
 \frac{q^{\frac1{16}}
(-qt;q)_\infty(-t^{-1};q)_\infty (q;q)_\infty^2}{2
 (qt;q)_\infty(t^{-1};q)_\infty (-q;q)_\infty}
 \end{aligned}
 $$
where we have used a version of the Jacobi triple product identity
$$\sum_{k\in\Z+\hf} q^{\frac{k^2}{2}} t^k =q^{\frac18}t^\hf
(q;q)_\infty (-qt;q)_\infty(-t^{-1};q)_\infty.$$
Comparing the two formulas of $1$-point function, we have the
following.

\begin{cor}
The following $q$-identity holds:
$$\begin{aligned}
  & \frac{ (-qt;q)_\infty(-t^{-1};q)_\infty (q;q)_\infty^2}{
 (qt;q)_\infty(t^{-1};q)_\infty (-q;q)_\infty^2}
  \\ &
\qquad =\frac{t+1}{t-1}
 + 2\sum_{r =0}^\infty
 \left[ \frac{(-1)^r q^{r+1} t}{1 -q^{r+1} t}
  -\frac{(-1)^r q^{r+1} t^{-1}}{1 -q^{r+1} t^{-1}}
 \right].
  \end{aligned}
$$
The right hand side of the identity is known (cf. \cite{W2}) to be
equal to
$$
2 t \frac{d}{dt}
 \ln
  \left(
   t^{-\hf}
\frac{(t; q^2)_\infty (q^2t^{-1}; q^2)_\infty}{(qt; q^2)_\infty
(qt^{-1} ; q^2)_\infty}
 \right).
 $$
\end{cor}

\subsection{The $n$-point $\binf$-functions of level $l+\hf$}

The character of the irreducible $Spin(2l+1)$-module associated to
$\la\in \Sigma(Pin)$ is also given by
$\text{ch}^b_\la(z_1,\dots,z_l)$ in (\ref{eq:charB}).

\begin{lemma}  \label{lem:howeBb}
We have the following $q$-series identity:
\begin{align*}
  \mathfrak B^{\hf}_{(\hf)} &(q,{\bf t}) \cdot \prod_{i=1}^l \bff_b
(z_i,q; t_1,\dots,t_n)
\\
&= \sum_{\la \in\Sigma(Pin)} \text{ch}^b_\la(z_1,\dots,z_l) \cdot
\mathfrak B_\la^{l+\half}(q;t_1,\ldots,t_n).
\end{align*}
\end{lemma}

\begin{proof}
Follows from the Howe duality in Proposition~\ref{howeBb} and
(\ref{eq:Bhf}). Note that the cancellation of a factor $2$ has
occurred.
\end{proof}

\begin{theorem} \label{th:Blplushalf}
The $n$-point $\binf$-correlation function of level $l+\hf$,
$\mathfrak B^{l+\hf}_\la(q,t_1,\ldots,t_n)$, is equal to
\begin{eqnarray*}
 && \mathfrak B^{\hf}_{(\hf)} (q;{\bf t}) \times \\
 &&\times
 \sum_{\sigma\in W(B_l)}\left(-1\right)^{\ell(\sigma)}
q^{\frac{\Vert \la+\rho-\sigma(\rho) \Vert^2}{2}} \prod_{a=1}^l
 \Big(
\sum_{\vec{\eps}_a\in\{\pm 1\}^n}
 [\vec{\eps}_a] (\Pi {\bf t}^{\vec{\eps}_a})^{k_a}
 F_{bo}(q;{\bf t}^{{\vec{\eps}_a}})
 \Big)
\end{eqnarray*}
where $k_a = (\la+\rho-\sigma(\rho), \ep_a)$.
\end{theorem}

\begin{proof}[Sketch of a proof]
The proof is completely parallel to the proof for
Theorem~\ref{th:lplushalf}, now with the help of (\ref{eq:charB}),
(\ref{eq:denomB}), Lemmas~\ref{weylcharB}, \ref{lem:nptB} and
\ref{lem:howeBb}.
\end{proof}
\subsection{The $q$-dimension of a $\binf$-module of level $l+\hf$}

In the same manner as in Section~\ref{sec:qdim-dinf-l}, we can
derive the following $q$-dimension formula from the Howe duality
in Proposition~\ref{howeBb}. The second formula below is obtained
from the first one by using Lemma~\ref{lem:weyldenom} and the
explicit root system of type $B_l$.
%
\begin{prop}
For $\la \in \Sigma(Pin)$, we have
\begin{align*}
  \dim_q & L(\binf;\Lambda'(\la))    \\
 & =
q^{\frac1{16}}\frac{(-q;q)_\infty}{(q;q)_\infty^l} \cdot
\sum_{\sigma\in W(B_l)}(-1)^{\ell(\sigma)} q^{\frac{\Vert
\la+\rho-\sigma(\rho) \Vert^2}{2}} \\
 &=
q^{\frac1{16}}\frac{(-q;q)_\infty}{(q;q)_\infty^l} \cdot
q^{\frac{\Vert\la\Vert^2}{2}}
  \prod_{1\leq i\leq l} \left(1-q^{\la_i+l-i+1/2}\right) \times\\
 & \qquad  \times \prod_{1\leq i<j\leq l}
\left(1-q^{\la_i-\la_j+j-i}\right)
\left(1-q^{\la_i+\la_j+2l-i-j+1}\right).
\end{align*}
\end{prop}

\end{document}